\documentclass[11pt]{amsart}

\usepackage{amsmath,amssymb}
\usepackage{indentfirst}
\usepackage[all]{xy}
\newtheorem{tm}{Theorem}[section]
\newtheorem{lemma}[tm]{Lemma}

\newtheorem{defi}[tm]{Definition}

\newtheorem{cor}[tm]{Corollary}

\newtheorem{remark}[tm]{Remark}

\begin{document}

\renewcommand{\P}{\mathbb{P}}
\newcommand{\R}{\mathbb{R}}
\newcommand{\N}{\mathbb{N}}
\newcommand{\Z}{\mathbb{Z}}
\newcommand{\Q}{\mathbb{Q}}
\newcommand{\scf}{\mathcal{F}}
\newcommand{\scc}{\mathcal{C}}
\newcommand{\scs}{\mathcal{S}}
\newcommand{\scu}{\mathcal{U}}
\newcommand{\scv}{\mathcal{V}}
\newcommand{\ra}{\rightarrow}
\newcommand{\ov}{\overline}
\newcommand{\bs}{\backslash}
\renewcommand{\mid}{\big\vert}
\newcommand{\al}{\alpha}
\newcommand{\da}{\delta}
\newcommand{\e}{\varepsilon}
\newcommand{\et}{\emptyset}
\newcommand{\ga}{\gamma}
\newcommand{\vp}{\varphi}
\newcommand{\Sa}{\Sigma}
\newcommand{\sa}{\sigma}
\newcommand{\la}{\lambda}
\renewcommand{\mod}{\mathrm{mod}}

\newcommand{\mesh}{\mathop{\mathrm{mesh}}}
\newcommand{\ord}{\mathop{\mathrm{ord}}}
\newcommand{\card}{\mathop{\mathrm{card}}}
\newcommand{\diam}{\mathop{\mathrm{diam}}}
\newcommand{\im}{\mathop{\mathrm{im}}}
\newcommand{\st}{\mathop{\mathrm{st}}}
\newcommand{\ANR}{\mathop{\mathrm{ANR}}}
\renewcommand{\int}{\mathop{\mathrm{int}}}
\newcommand{\inv}{^{-1}}
\newcommand{\Sh}{\mathop{\mathrm{Sh}}}
\newcommand{\UV}{\mathop{\mathrm{UV}}}
\newcommand {\Tor}{\mathop{\mathrm{Tor}}}

\title{Bockstein basis and resolution theorems in\\
extension theory\\ \footnotesize{\textnormal{published in \emph{Topology and its Applications} 157 (2010)  674-691}}}
\author{Vera Toni\'{c}}
\address{Department of Mathematics\\
University of Oklahoma\\
601 Elm Ave, room 423\\
Norman, Oklahoma 73019\\
USA}\email{vtonic@ou.edu}

\date{22 November 2009}


\begin{abstract}
We prove a generalization of the Edwards-Walsh Resolution Theorem:

\noindent \textbf{Theorem}: \emph{Let $G$ be an abelian group with
$P_G=\mathbb{P}$, where $P_G=\{ p \in \P: \Z_{(p)}\in$ Bockstein Basis $ \sigma(G)\}$. Let $n\in \N$ and let $K$ be a connected \emph{CW}-complex with $\pi_n(K)\cong G$, $\pi_k(K)\cong 0$ for $0\leq k< n$. Then for every compact metrizable
space $X$ with $X\tau K$ (i.e., with $K$ an absolute extensor for
$X$), there exists a compact metrizable space $Z$ and a surjective
map $\pi: Z \rightarrow X$ such that
\begin{enumerate}
\item[(a)] $\pi$ is cell-like,
\item[(b)] $\dim Z\leq n$, and
\item[(c)] $Z\tau K$.
\end{enumerate}
}
\end{abstract}


\maketitle \markboth{V.\ Toni\'{c}}{Bockstein Basis and Resolution Theorems in Extension Theory}

\section{Introduction}

The objective of this paper will be to prove the following resolution theorem:
\begin{tm}\label{T} Let $G$ be an abelian group with
$P_G=\mathbb{P}$, where $P_G=\{ p \in \P: \Z_{(p)}\in$ Bockstein Basis $ \sigma(G)\}$. Let $n\in \N$ and let $K$ be a connected \emph{CW}-complex with $\pi_n(K)\cong G$, $\pi_k(K)\cong 0$ for $0\leq k< n$. Then for every compact metrizable
space $X$ with $X\tau K$ (i.e., with $K$ an absolute extensor for
$X$), there exists a compact metrizable space $Z$ and a surjective
map $\pi: Z \rightarrow X$ such that
\begin{enumerate}
\item[(a)] $\pi$ is cell-like,
\item[(b)] $\dim Z\leq n$, and
\item[(c)] $Z\tau K$.
\end{enumerate}
\end{tm}

The word resolution refers to a map between topological spaces where the domain is in some way better than the range, and the fibers (point preimages) meet certain requirements.

Let us look at some examples of resolution theorems. Here is the cell-like resolution theorem, first stated by R.\ Edwards (\cite{Ed}), and later proven by J.\ Walsh in \cite{Wa}:

\begin{tm}\emph{(R.~Edwards - J.~Walsh, 1981) \cite{Wa}}:\label{EdWa}
For every compact metrizable space $X$ with $\dim_{\Z} X
\leq n$, there exists a compact metrizable space $Z$ and a
surjective map $\pi :Z \ra X$ such that
$\pi$ is cell-like, and $\dim Z \leq n$. \hfill $\square$
\end{tm}

If $n\in \N$, then a \emph{subset} $Y \subset \R^n$ is called \emph{cellular} if $Y$ can be written as the intersection of a nested collection of $n$-cells in $\R^n$. A \emph{space} $Y$ is called \emph{cell-like} if for some $n\in \N$, there is an embedding $F:Y \ra \R^n$ so that $F(Y)$ is cellular.
A \emph{map} $\pi: Z \ra X$ is called \emph{cell-like} if for each $x\in X$, $\pi^{-1}(x)$ is cell-like. Whenever $X$ is a finite-dimensional compact metrizable space, then $X$ is cell-like if and only if $X$ has the shape of a point. To detect that a compact metrizable space has the shape of a point, it is sufficient to prove that there is an inverse sequence $(Z_i,p_i^{i+1})$ of compact metrizable spaces $Z_i$ whose limit is homeomorphic to $X$ and such that for each $i \in \N$, $p_i^{i+1}:Z_{i+1} \ra Z_i$ is null-homotopic. It is also sufficient to show that every map of $X$ to a CW-complex is null-homotopic.

The Edwards-Walsh Theorem has been generalized to the class of arbitrary metrizable spaces by L.\ Rubin and P.\ Schapiro (\cite{RS1}), and to the class of arbitrary compact Hausdorff spaces by S.\ Marde\v{s}i\'{c} and L.\ Rubin (\cite{MR}).

A similar statement to the Edwards-Walsh Theorem was proven by A.\ Dranishnikov, for the group $\Z/p$, where $p$ is an arbitrary prime number:

\begin{tm} \label{Dran} \emph{(A.\ Dranishnikov, 1988) \cite{Dr2}}: For
every compact metrizable space $X$ with $\dim_{\Z/p} X$ $\leq n$,
there exists a compact metrizable space $Z$ and a surjective map
$\pi :Z \ra X$ such that $\pi$ is $\Z/p$-acyclic, and $\dim Z \leq n$. \hfill $\square$
\end{tm}

A map $\pi :Z \ra X$ between topological spaces is called $G$-\emph{acyclic} if
all its fibers $\pi^{-1}(x)$ have trivial reduced \v{C}ech
cohomology with respect to the group $G$, or, equivalently,  every map
$f: \pi^{-1}(x) \ra K(G,n)$ is nullhomotopic.
Note that a map $\pi :Z \ra X$ being cell-like implies that $\pi$ is also $G$-acyclic.

\vspace{4mm}

Akira Koyama and Katsuya Yokoi (\cite{KY1}) were able to obtain this $\Z/p$-resolution theorem of Dranishnikov both for the class of metrizable spaces and for the class of compact Hausdorff spaces. Dranishnikov proved a statement similar to Theorem~\ref{Dran} for the group $\Q$ (\cite{Dr4}), but he could only obtain $\dim Z \leq n+1$, and if $n\geq 2$, then additionally $\dim_\Q Z \leq n$. This result was later improved by M.\ Levin:

\begin{tm} \label{Lev} \emph{(M.\ Levin, 2005) \cite{Le2}}: Let $n\in \N_{\geq2}$. Then for
every compact metrizable space $X$ with $\dim_{\Q} X$ $\leq n$,
there exists a compact metrizable space $Z$ and a surjective map
$\pi :Z \ra X$ such that $\pi$ is $\Q$-acyclic, and $\dim Z \leq n$. \hfill $\square$
\end{tm}

The obvious question was whether a theorem similar to Theorem~\ref{Dran}
could be stated for compact metrizable spaces and arbitrary abelian groups.
In their work \cite{KY2}, Koyama and Yokoi made a substantial amount of progress in answering this question. Their method relied heavily on the existence of Edwards-Walsh complexes, which have been studied by J.\ Dydak and J.\ Walsh in \cite{DW}, and which had been applied originally, in a rudimentary form, in \cite{Wa}.
However, using a different approach from the one in \cite{KY2}, M.\ Levin has proved a very strong generalization for Theorems \ref{EdWa} and  \ref{Dran}, concerning compact metrizable spaces and arbitrary abelian groups:

\begin{tm}\label{Lev2} \emph{(M.\ Levin, 2003) \cite{Le1}}: Let
$G$ be an abelian group and let $n \in \N_{\geq 2}$.
Then for every compact metrizable space $X$ with $\dim_G X \leq
n$, there exists a compact metrizable space $Z$ and a surjective
map $\pi :Z \rightarrow X$ such that:

\emph{(a)} $\pi$ is $G$-acyclic,

\emph{(b)} $\dim Z \leq n+1$, and

\emph{(c)} $\dim_G Z \leq n$. \hfill $\square$
\end{tm}

The requirement of $n \in \N_{\geq 2}$ in Levin's Theorem cannot be improved because there is a counterexample for $n=1$ ($G=\Q$, \cite{Le1}). The requirement that $\dim Z \leq n+1$ cannot be improved either -- there is a counterexample for $\dim Z \leq n$ ($G=\Z/p^{\infty}$, \cite{KY2}).
The part that may be improved is $\dim_G X \leq n$, using the characterization of cohomological dimension by extension of maps.
Namely, for any paracompact Hausdorff space $X$, any abelian group $G$ and $n\in \N$, $\dim_G X \leq n$ if and only if every map of a closed subspace of $X$ to   $K(G,n)$ can be extended to a map of $X$ to $K(G,n)$. By $K(G,n)$ we will always mean an Eilenberg-MacLane CW-complex of type $(G,n)$, and such is characterized (up to homotopy equivalence) by having $\pi_n\cong G$ and  $\pi_k$ trivial for all other $k$.

This fact about extending maps from any closed subspace of $X$ to a $K(G,n)$ can be written as $K(G,n) \in AE(X)$ ($K(G,n)$ is an \emph {absolute extensor} for $X$). Another notation, and the one we will be using, is $X\tau K(G,n)$. In fact, for any two topological spaces $X$ and $Y$, $X\tau Y$ will mean that every map from a closed subspace of $X$ to $Y$ can be extended continuously over $X$.

\vspace{4mm}

So, in order to generalize the requirement  $\dim_G X \leq n$ from  Theorem~\ref{Lev2}, note that $\dim_G X \leq n$
$\Leftrightarrow$ $X \tau K(G,n)$, and  replace a $K(G,n)$ with a CW-complex upon which the demands will be less strict.
Here is a theorem generalizing  Theorem~\ref{Lev2} for some abelian groups.

\begin{tm} \label{RuSc}\emph{(L.\ Rubin - P.\ Schapiro, 2005) \cite{RS2}}:
Let $G$ be an abelian group with $P_G\neq\P$, where
$P_G=\{ p \in \P: \Z_{(p)} \in$ Bockstein basis $ \sa(G)\}$. Let
$n\in \N_{\geq 2}$, and let $K$ be a connected \emph{CW}-complex with
$\pi_n(K)\cong G$, $\pi_k(K)\cong 0$ for $0\leq k< n$. Then for every
compact metrizable space $X$ with $X\tau K$, there exists a
compact metrizable space $Z$ and a surjective map $\pi: Z
\rightarrow X$ such that:

\emph{(a)} $\pi$ is $G$-acyclic,

\emph{(b)} $\dim Z \leq n+1$, and

\emph{(c)} $Z\tau K$.\hfill $\square$


\end{tm}

Note that  the statement of Theorem~\ref{RuSc} does not cover the case when $P_G=\P$. In fact, the statement of this theorem will be true when $P_G=\P$, but in this case the statement can be improved, as shown in Theorem~\ref{T}.

The author wishes to thank Dr.\ Leonard Rubin and the referee of this paper, for their wise suggestions for improvements in all upcoming sections.

\vspace{1mm}

Before we proceed,  let us review some basic facts from Bockstein theory.

\section{Bockstein Theory}

The cohomological dimension of a given compact metrizable space depends on the coefficient group, which can be any abelian group
and there are uncountably many of them. It turns out that in the case of compact metrizable spaces,
it suffices to consider only countably many groups.
M.\ F.~Bockstein found  an algorithm for
computation of the cohomological dimension with respect to a given abelian group $G$ by
means of cohomological dimensions with coefficients taken from a countable family of abelian groups $\sa(G)$.
His definition of $\sa(G)$ was also used by V.\ I.~Kuz'minov (\cite{Ku}), and later adapted by E.~Dyer (\cite{Dy}), and then by A.~Dranishnikov (\cite{Dr3}).

Thus there are three different definitions of a Bockstein basis $\sa(G)$, which are not equivalent in general, but which are equivalent from the point of view of cohomological dimension. This can be shown using the Bockstein Theorem and Bockstein Inequalities, which will be stated in this section.

\vspace{4mm}

\noindent \textbf{Notation}:
\begin{enumerate}
\item[(1)] $\P$ stands for the set of all prime numbers,
\item[(2)] $\Z_{(p)}=\{ \frac{m}{n} \in \Q \
: \ n$ is not divisible by $p \}$ is called the $p$-\emph{localization of the integers}, and
\item[(3)] $\Z/p^\infty =\{ \frac{m}{n} \in \Q/\Z \ : \ n=p^k$ for some $k\geq 0\}$ is called \emph{the quasi-cyclic $p$-group}.
\end{enumerate}

\vspace{2mm}

For an abelian group $G$, we say that  \emph{an element $g\in G$ is divisible by} $n \in \Z \setminus\{0\}$ if the equation $nx=g$ has a solution in $G$, \emph{$G$ is divisible by $n$} if all of its elements are divisible by $n$, and  $G$ is a \emph{divisible group} if $G$ is divisible by all $n \in \Z \setminus \{0\}$.

For an abelian group $G$, $\Tor G$ is the subgroup of all elements of $G$ of finite order, and $p$--$\Tor G$ is the subgroup of all elements whose order is a power of $p$, that is, $p$--$\Tor G = \{g \in G : p^kg=0 \ \mathrm{for} \ \mathrm{some} \ k\geq 1\}$.

Here is the definition of a Bockstein basis $\sa(G)$ that we will use,
adapted from the original one by E.~Dyer (\cite{Dy}).

\begin{defi} Let $G$ be an abelian
group, $G\neq 0$. Then $\sa (G)$ is the subset of $\{ \Q \} \cup \{
\Z/p,\Z/p^\infty,\Z_{(p)} : p \in \mathbb{P}\}$ defined by:

\vspace{2mm}
\begin{tabular}{rccl}
\hspace*{-1cm} \emph{(\textbf{I})} & $\Q \in \sa(G)$ \hspace{5mm} & $\Leftrightarrow$ & $G$ contains an element of infinite order\\
                            &                 & $\Leftrightarrow$ & $G/\Tor G \neq 0$\\
\hspace*{-1cm} \emph{(\textbf{II})} &  $\Z_{(p)} \in \sa(G)$ & $\Leftrightarrow$ &
$G$ satisfies the following: $\exists g \in G$ such that $\forall k\in \Z_{\geq0}$,\\
                            &                &               &  $p^kg$ is not divisible by $p^{k+1}$\\
                            &                 & $\Leftrightarrow$ & $G/\Tor G$ is not divisible by $p$\\
\end{tabular}

\begin{tabular}{rccl}
\hspace*{-1cm} \emph{(\textbf{III})} &  $\Z/p \in \sa(G)$  & $\Leftrightarrow$ &   $G$
contains an element of order $p^k$, for some $k\in \N$,\\
                            &                &               & which is not divisible by $p$\\
                            &                 & $\Leftrightarrow$ & $p$--$\Tor G$ is not divisible by $p$\\
\hspace*{-1cm} \emph{(\textbf{IV})} & $\Z/p^\infty \in \sa(G)$  & $\Leftrightarrow$ &
                                             $p$--$\Tor G \neq 0$ and $p$--$\Tor G$ is divisible by $p$.

\end{tabular}
\end{defi}

\vspace{2mm}

\begin{tm}[Bockstein Inequalities] \emph{\cite{Dr3}}:
For any compact metrizable space X the following inequalities hold:
\begin{enumerate}
\item[(BI1)] $\dim_{\Z/p^\infty} X \leq \dim_{\Z/p} X$,
\item[(BI2)] $\dim_{\Z/p} X \leq \dim_{\Z/p^\infty} X +1$,
\item[(BI3)] $\dim_{\Z/p} X \leq \dim_{\Z_{(p)}} X$,
\item[(BI4)] $\dim_{\Q} X \leq \dim_{\Z_{(p)}} X$,
\item[(BI5)] $\dim_{\Z_{(p)}} X \leq \max \{\dim_{\Q} X ,\dim_{\Z/p^\infty} X +1\}$,
\item[(BI6)] $\dim_{\Z/p^\infty} X \leq \max \{\dim_{\Q} X ,\dim_{\Z_{(p)}} X -1\}.\hfill \square$
\end{enumerate}
\end{tm}

\begin{tm}[Bockstein Theorem] \emph{\cite{Dy}}: If $G$ is an abelian group and $X$ is a locally compact space, then \
$\displaystyle{\dim_G X = \sup_{H\in \sa(G)} \dim_H X.}$\ \hfill $\square$
\end{tm}

\vspace{2mm}

Now let $P_G:=\{ p \in \P: \Z_{(p)} \in \sa(G)\}$.

\begin{lemma}\label{PGPrimesZ}
If $G$ is an abelian group such that $P_G=\P$, then for any compact metrizable space $X$, $\dim_G X = \dim_\Z X$.
\end{lemma}
\noindent \textit{Proof}: $P_G=\P$ means that for each $p \in \P$, $\Z_{(p)} \in\sa(G)$. By the Bockstein Inequalities (BI4), (BI3) and (BI1), the supremum
$\displaystyle{\sup_{H\in \sa(G)} \dim_H X}$ has to be achieved at  $\displaystyle{\sup_{p\in\P} \ \dim_{\Z_{(p)}}} X$. Since $\sa(\Z)=\{\Q\} \cup \{\Z_{(p)} : p \in \P\}$,\ \
$\displaystyle{\sup_{H\in \sa(G)} \dim_H X= \sup_{H\in \sa(\Z)} \dim_H X}$. \hfill $\square$

\section{Walsh Technical Lemma and Edwards Type Theorem}

\noindent This will be a statement needed to produce a resolution $\pi : Z \twoheadrightarrow X$, based on \cite{Wa}.

\noindent\textbf{Notation}: $B_r(x)$ stands for the \textit{closed} ball with radius $r$, centered at $x$.

\begin{lemma}[\textbf{Generalized Walsh Lemma}]\label{genWa} Let $\mathbf{X}
=(P_i,f^{i+1}_i)$ be an inverse sequence of compact metric spaces $(P_i,d_i)$ of diameter less than $1$,
$\mathbf{Z}=(M_i,g_i^{i+1})$  an inverse sequence of Hausdorff compacta,
$X= \lim \mathbf{X}$ and $Z=\lim \mathbf{Z}$. Assume also that we have maps $\phi_i:M_i \ra P_i$, and,
for each $i \in \N$ we have numbers
$0 < \e(i)<\frac{\da(i)}{3}<1$, satisfying:
\begin{enumerate}
\item[(I)] for $i\geq 2$, $\phi_{i-1}\circ g_{i-1}^{i}$ and
$f_{i-1}^{i}\circ \phi_i$ are $\frac{\e(i-1)}{3}$ - close,

\item[(II)] for $i\geq 2$ and for any $y \in P_{i}$, $\diam
\ (f_{i-1}^{i}(B_{\da(i)}(y))) \ < \ \frac{\e(i-1)}{3}$, and

\item[(III)] for $i>j$ and for any $y \in P_{i}$, $\ \diam \
(f_j^i(B_{\e(i)}(y))) \ < \ \frac{\e(j)}{2^i}$.

\end{enumerate}
Then there is a map $\pi: Z \ra X$ such that for all $x = (x_i)\in X$:
\begin{enumerate}
\item[(IV)]
$\pi^{-1}(x)\ = \ \lim \ (\phi_i^{-1}(B_{\da(i)}(x_i)) ,
g_i^{i+1})\ =\ \lim \ (\phi_i^{-1}(B_{\e(i)}(x_i)) , g_i^{i+1})$
\end{enumerate}
(here $g_i^{i+1}$ stands for the appropriate restriction).

\vspace{2mm}

\noindent If, in addition, we have that:
\begin{enumerate}
\item [(V)] for all $x=(x_i) \in X$ and for all $i$, $\phi_i^{-1}(B_{\e(i)}(x_i))\neq \et$,
\end{enumerate}
then $\pi^{-1}(x) \neq \et$, so the map
$\pi$ will be surjective.
\end{lemma}

\noindent \textit{Proof}: The following diagram will help in visualizing the steps of this proof.
\begin{displaymath}
\xymatrix{
\cdots & M_i\ar[l] \ar[d]_{\phi_i} &  M_{i+1} \ar[l]_{g_i^{i+1}} \ar[d]^{\phi_{i+1}} & \cdots\ar[l] & Z\ar@{-->}[d]^\pi\\
\cdots & P_i\ar[l] &  P_{i+1} \ar[l]^{f_i^{i+1}} & \cdots\ar[l]   & X
}
\end{displaymath}

\vspace{2mm}

 Let $z=(z_i)$ be an element of $Z\subset \prod_{i=1}^\infty
M_i$; so $g_i^{i+1}(z_{i+1})=z_i$ and $\phi_i(z_i) \in P_i$, for all $i\in \N$. Define a sequence in
$\prod_{i=1}^\infty P_i$ as follows:
\begin{align*}
x^1 &= (\phi_1(z_1),\phi_2(z_2),\phi_3(z_3),\phi_4(z_4),\ldots)\\
x^2 &= (f_1^2(\phi_2(z_2)),\phi_2(z_2),\phi_3(z_3),\phi_4(z_4),\ldots)\\
x^3 &= (f_1^3(\phi_3(z_3)),f_2^3(\phi_3(z_3)),\phi_3(z_3),\phi_4(z_4),\ldots)\\
\vdots\\
x^j &= (f_1^j(\phi_j(z_j)),f_2^j(\phi_j(z_j)),\ldots,f_{j-1}^j(\phi_j(z_j)),\phi_j(z_j),\phi_{j+1}(z_{j+1}),\ldots)\\
x^{j+1} &= (f_1^{j+1}(\phi_{j+1}(z_{j+1})),f_2^{j+1}(\phi_{j+1}(z_{j+1})),\ldots,f_{j}^{j+1}(\phi_{j+1}(z_{j+1})),\phi_{j+1}(z_{j+1}),\phi_{j+2}(z_{j+2}),\ldots)\\
\vdots
\end{align*}

Let $\pi_j:Z \ra \prod_{i=1}^\infty P_i$ be defined by
$\pi_j(z):=x^j$. Note that $\pi_j$ are continuous because coordinate maps $x^j$ are continuous. We shall employ the metric $d$ on $\prod_{i=1}^\infty P_i$ given by
$$d((s_i),(r_i)):=\sum_{i=1}^{\infty}\frac{d_i(s_i,r_i)}{2^i}.$$ We would like to show that
$(\pi_j(z))_{_{j\in\N}}$ is a Cauchy sequence in
$\prod_{i=1}^\infty P_i$. Properties we will need are:
\begin{enumerate}
\item[(1)] for $j\geq 2$, $f_{j-1}^j(\phi_{j}(z_j))$ and
$\phi_{j-1}(z_{j-1})=\phi_{j-1}(g_{j-1}^j(z_j))$ are $\e(j-1)$-close, and

\item[(2)] for $i>j$, $f_j^{i+1}(\phi_{i+1}(z_{i+1}))$ and $f_j^i(\phi_i(z_i))$ are
$\frac{\e(j)}{2^i}$-close.
\end{enumerate}

Property (1) follows from (I). Property (2) is true because: by
(1)$_{i+1}$, $f_{i}^{i+1}(\phi_{i+1}(z_{i+1}))$ and $\phi_i(z_i)$ are $\e(i)$-close,
so $f_{i}^{i+1}(\phi_{i+1}(z_{i+1})) \in B_{\e(i)}(\phi_i(z_i))$. Therefore
$f_{j}^{i+1}(\phi_{i+1}(z_{i+1}))=$ $f_j^i(f_i^{i+1}(\phi_{i+1}(z_{i+1}))) \in
f_j^i(B_{\e(i)}(\phi_i(z_i)))$, and $\diam
f_j^i(B_{\e(i)}(\phi_i(z_i)))<\frac{\e(j)}{2^i}$, by (III). So
$f_j^{i+1}(\phi_{i+1}(z_{i+1}))$ and $f_j^i(\phi_i(z_i))$ are
$\frac{\e(j)}{2^i}$-close.

Note that by
(2)$_{j>q}$ and (1)$_{j+1}$,
\begin{align*}
d(\pi_j(z),\pi_{j+1}(z)) &= \left(\sum_{q=1}^{j-1}
\frac{d_q(f_q^j(\phi_j(z_j)),f_q^{j+1}(\phi_{j+1}(z_{j+1})))}{2^q} \right)+
\frac{d_j(\phi_j(z_j),f_j^{j+1}(\phi_{j+1}(z_{j+1})))}{2^j}\ \\
 &<\left(\sum_{q=1}^{j-1}\frac{\e(q)}{2^j}\frac{1}{2^q}\right)\ +\
 \frac{\e(j)}{2^j} \ <\
 \frac{1}{2^j}\left(\sum_{q=1}^{j-1}\frac{1}{2^q}\right)\ +\
 \frac{1}{2^j}\\
 &< \frac{1}{2^j}\left(\left(\sum_{q=1}^{\infty}\frac{1}{2^q}\right)\ +\
 1\right)\ =\ \frac{1}{2^{j-1}}.
\end{align*}
Therefore, for the indexes $j$ and $j+k$ we get:
\begin{align*}
d(\pi_j(z),\pi_{j+k}(z)) &\leq d(\pi_j(z),\pi_{j+1}(z)) +
d(\pi_{j+1}(z),\pi_{j+2}(z)) + \ldots +
d(\pi_{j+k-1}(z),\pi_{j+k}(z))\\
&<\frac{1}{2^{j-1}}+\frac{1}{2^{j}}+\ldots+\frac{1}{2^{j+k-2}}<\frac{1}{2^{j-2}}\cdot\sum_{i=1}^{\infty}\frac{1}{2^i}\
=\frac{1}{2^{j-2}}.
\end{align*}
Thus $(\pi_j(z))_{_{j\in\N}}$ is a Cauchy sequence in the compact
metric space $\prod_{i=1}^\infty P_i$, and therefore it is
convergent. Define $\pi(z):= \displaystyle{\lim_{j\ra\infty}}
\pi_j(z)$.

Notice that for any $k \in \N$, and for any $z\in Z$,
$$d(\pi_k(z),\pi(z)) \leq \sum_{j=k}^\infty
d(\pi_j(z),\pi_{j+1}(z)) < \sum_{j=k}^\infty\frac{1}{2^{j-1}}
=\frac{1}{2^{k-2}}.$$ So the sequence $(\pi_j)_{j\in\N}$ converges
uniformly to $\pi$. Therefore $\pi: Z \ra \prod_{i=1}^\infty P_i$ is a
continuous function.

\vspace{2mm}

We would like to see that $\pi(Z)\subset X$. If $y_j$ is $j$-th
coordinate of $\pi(z)$ for some $z\in Z$, then $y_j =
\displaystyle{\lim_{i>j}}\ f_j^i(\phi_i(z_i))$. Therefore if $j>1$,
\begin{align*}
f_{j-1}^j(y_j) & =f_{j-1}^j(\lim_{i>j}\
f_j^i(\phi_i(z_i)))=\lim_{i>j}(f_{j-1}^j(
f_j^i(\phi_i(z_i))))\\
& = \lim_{i>j}(f_{j-1}^i(\phi_i(z_i)))=\lim_{i>j-1}(f_{j-1}^i(\phi_i(z_i)))=y_{j-1}.
\end{align*}
So $\pi(z)\in X$, i.e., $\pi(Z)\subset X$.

Now that we have a map $\pi:Z\ra X$, we need to see what its
fibers are. Take any $x=(x_i)\in X$. From (II)$_i$ and (I)$_i$, we
will get that
\begin{enumerate}
\item[(3)] $g_{i-1}^i(\phi_i^{-1}(B_{\da(i)}(x_i))) \ \subset \
\phi_{i-1}^{-1}(B_{\e(i-1)}(x_{i-1}))$.
\end{enumerate}
Here is why: take any $y \in \phi_i^{-1}(B_{\da(i)}(x_i))$, i.e., $\phi_i(y) \in B_{\da(i)}(x_i)$.
Note that (II)$_i$:\\ $\diam \ (f_{i-1}^{i}(B_{\da(i)}(x_i)))  <
\frac{\e(i-1)}{3}$.  Hence $d_{i-1}(f_{i-1}^i(\phi_i(y)),
f_{i-1}^i(x_i))<\frac{\e(i-1)}{3}$,\linebreak i.e.,
$d_{i-1}(f_{i-1}^i(\phi_i(y)),x_{i-1})<\frac{\e(i-1)}{3}$. By (I)$_i$:
$d_{i-1}(\phi_{i-1}(g_{i-1}^i(y)), f_{i-1}^i(\phi_i(y)))<\frac{\e(i-1)}{3}$, and
therefore
\begin{align*}
d_{i-1}(x_{i-1}, \phi_{i-1}(g_{i-1}^i(y)))& \leq d_{i-1}(x_{i-1},
f_{i-1}^i(\phi_i(y))) + d_{i-1}(f_{i-1}^i(\phi_i(y)),
\phi_{i-1}(g_{i-1}^i(y)))\\
& <\frac{2\e(i-1)}{3}<\e(i-1).
\end{align*}
So $\phi_{i-1}(g_{i-1}^i(y))\in B_{\e(i-1)}(x_{i-1})$, and therefore $g_{i-1}^i(y) \in \phi_{i-1}^{-1}( B_{\e(i-1)}(x_{i-1}))$, so (3) is true.

\vspace{2mm}

As a consequence of (3) and the fact that $\e(i)<\da(i)$, both\\ $(\phi_i^{-1}(B_{\da(i)}(x_i)),
g_{i-1}^{i}\vert_{\phi_i^{-1}(B_{\da(i)}(x_i))})$ and
$(\phi_i^{-1}(B_{\e(i)}(x_i)) ,
g_{i-1}^{i}\vert_{\phi_i^{-1}(B_{\e(i)}(x_i)) })$ are inverse
sequences with the same limit. Now we would like to show that this
limit is $\pi^{-1}(x)$.

Let us show that $\lim (\phi_i^{-1}(B_{\e(i)}(x_i)) ,
g_{i-1}^{i}) \subset  \pi^{-1}(x)$, where $g_{i-1}^i$ stands for
the appropriate restriction. Take any $z=(z_i) \in \lim (\phi_i^{-1}(B_{\e(i)}(x_i)) ,
g_{i-1}^{i})$. Note that
\begin{enumerate}
\item[(4)] the $j$-th
coordinate of $\pi(z)$ is $\displaystyle{\lim_{i>j}}\ f_j^i(\phi_i(z_i))$.
\end{enumerate}
Since $z_i\in \phi_i^{-1}(B_{\e(i)}(x_i))$, we have that $\phi_i(z_i) \in B_{\e(i)}(x_i)$.
Condition (III)$_i$, which says that $ \diam \
(f_j^i(B_{\e(i)}(x_i)))$ $<  \frac{\e(j)}{2^i}$, implies that
$f_j^i(\phi_i(z_i))$ and $x_j=f_j^i(x_i)$ are $\frac{\e(j)}{2^i}$-close.
Therefore $\displaystyle{\lim_{i>j}}\ f_j^i(\phi_i(z_i))=x_j$, so
$\pi(z)=x$, i.e., $z\in \pi^{-1}(x)$.

\vspace{2mm}

Let us demonstrate that $\pi^{-1}(x)\subset \lim (\phi_i^{-1}(B_{\da(i)}(x_i)) , g_{i-1}^{i})$. Suppose that $z=(z_i)\in Z$, and $z
\notin \lim (\phi_i^{-1}(B_{\da(i)}(x_i)) , g_{i-1}^{i})$. We
will show that $\pi(z)\neq x$.

Now $z \notin \lim (\phi_i^{-1}(B_{\da(i)}(x_i)) , g_{i-1}^{i})$
means that there is an index $j\in\N$ such that $z_j \notin
\phi_j^{-1}(B_{\da(j)}(x_j))$. So $d_j(\phi_j(z_j),x_j)>\da(j)$. The
inequality $\e(j)<\frac{\da(j)}{3}$ assures that
$B_{2\e(j)}(\phi_j(z_j))\cap B_{\e(j)}(x_j)=\emptyset$. If we look at the
distance between $\phi_j(z_j)$ and the $j$-th coordinate of $\pi(z)$ (see (4)), from
(1)$_{j+1}$ and (2)$_{k>j}$ we get:
\begin{align*}
d_j(\phi_j(z_j),\lim_{i>j}\ f_j^i(\phi_i(z_i)))\ &\leq
d_j(\phi_j(z_j),f_j^{j+1}(\phi_{j+1}(z_{j+1}))) \\
& \ \ +\sum_{k=j+1}^\infty d_j(f_j^{k}(\phi_k(z_{k})),f_j^{k+1}(\phi_{k+1}(z_{k+1}))) \\
&< \e(j)+\sum_{k=j+1}^\infty \frac{\e(j)}{2^k}\ =\
\e(j)+\frac{\e(j)}{2^j}\cdot\sum_{k=1}^\infty \frac{1}{2^k}\ <\
2\e(j).
\end{align*}
That is, the $j$-th coordinate of $\pi(z)$ is contained in
$B_{2\e(j)}(\phi_j(z_j))$, implying $\pi(z)\neq x$, i.e., $z\notin
\pi^{-1}(x)$.

So we get that
$$\lim (\phi_i^{-1}(B_{\e(i)}(x_i)), g_{i-1}^{i})
\subset  \pi^{-1}(x) \subset \lim (\phi_i^{-1}(B_{\da(i)}(x_i)),
g_{i-1}^{i}),$$
and since the left and right side of this
statement are equal, then (IV) is true.

If (V) is also true, i.e., $\pi^{-1}(x)$ is the inverse limit of an inverse sequence
of compact nonempty spaces,
then, according to  Theorem~2.4 from Appendix II of
\cite{Du},
$\pi^{-1}(x)\neq \et$. Thus, the map $\pi : Z \ra X$ is surjective. \hfill $\square$


\begin{remark}\label{unstable}
\emph{In some of the proofs that follow we will use stability theory, about which more details can be found in \S VI.1 of \cite{HW}. Namely, we will use the consequences of the Theorem VI.1. from \cite{HW}: if $X$ is a separable metrizable space with $\dim X \leq n$,  then for any map $f: X\ra I^{n+1}$ all values of $f$ are unstable. A point $y\in f(X)$ is called an  \emph{unstable value} of $f$ if for every $\da > 0$ there exists a map $g:X\ra I^{n+1}$ such that:
\begin{enumerate}
\item $d (f(x),g(x))<\da$ \ for every $x\in X$, and
\item$g(X)\subset I^{n+1}\setminus\{y\}$.
\end{enumerate}
Moreover, this map $g$ can be chosen so that $g=f$ on the complement of an arbitrary open neighborhood of $y$, and so that $g$ is homotopic to $f$ (see Corollary I.3.2.1 of \cite{MS1}).
}
\end{remark}

\begin{lemma}[\textbf{Special version of Walsh Lemma}] \label{specgenWa} Let $\mathbf{X}
=(P_i,f^{i+1}_i)$ be an inverse sequence of compact metric polyhedra $(P_i,d_i)$ with diameter less than $1$,
and let $L_i$ be triangulations of $P_i$.
Suppose that we have maps $g_i^{i+1}:\vert L_{i+1}^{(n+1)}\vert
\ra \vert L_{i}^{(n+1)}\vert$ such that $g_i^{i+1}(\vert L_{i+1}^{(n)}\vert) \subset \vert L_{i}^{(n)}\vert$, and let $\mathbf{Z}=(\vert
L_i^{(n)}\vert,g_i^{i+1})$ be the inverse sequence of subpolyhedra
$\vert L_i^{(n)}\vert \subset P_i$, where each $g_i^{i+1}$ stands
for the appropriate restriction. Let $X= \lim \mathbf{X}$, $Z=
\lim \mathbf{Z}$. Assume that for each $i \in \N$ we have numbers
$0 < \e(i)<\frac{\da(i)}{3}<1$, satisfying:
\begin{enumerate}
\item[(I)] for $i\geq 2$, $g_{i-1}^{i}\mid_{ \vert L_{i}^{(n)}\vert}$ and
$f_{i-1}^{i}\mid_{ \vert L_{i}^{(n)}\vert}$ are
$\frac{\e(i-1)}{3}$ - close,
%
\end{enumerate}
and conditions \emph{(II)} and \emph{(III)} from Lemma~\ref{genWa}.

Then there is a map $\pi: Z \ra X$ such that for all $x = (x_i)\in X$:
$$\pi^{-1}(x)\ =\ \lim \ (B_{\da(i)}(x_i)\cap \vert L_i^{(n)}\vert ,
g_i^{i+1})\ =\ \lim \ (B_{\e(i)}(x_i)\cap \vert L_i^{(n)} \vert, g_i^{i+1})$$
(here $g_i^{i+1}$ stands for the appropriate restriction).

\vspace{2mm}

\noindent If, in addition, we have that:
\begin{enumerate}
\item [(IV)] $\mesh L_i\ < \ \e(i)\ $, for all $i$,
\end{enumerate}
then for all $x \in X$ we have $\pi^{-1}(x) \neq \et$, so the map
$\pi$ will be surjective.

\noindent If we also have
\begin{enumerate}
\item [(V)] for $i\geq 1$ and for any $y \in P_i$,
$B_{\e(i)}(y)\subset P_{y,i}\subset B_{\da(i)}(y)$, where
$P_{y,i}$ is a contractible subpolyhedron of $\vert L_i \vert$,
and

\item[(VI)] for $i\geq 2$, $g_{i-1}^i(\vert
L_i^{(n+1)}\vert)\subset \vert L_{i-1}^{(n)}\vert$,
\end{enumerate}
then the map $\pi$ is cell-like.
\end{lemma}

\noindent \textit{Proof}: The following diagram will be useful.
\begin{displaymath}
\xymatrix{
\cdots & \vert L_i^{(n)}\vert\ar[l] \ar@{_{(}->}[d] & & \vert L_{i+1}^{(n)} \vert \ar[ll]_{g_i^{i+1}\vert_{\vert L_{i+1}^{(n)} \vert}} \ar@{_{(}->}[d] & \cdots\ar[l] & Z\ar@{-->}[d]^\pi\\
\cdots & P_i=\vert L_i\vert\ar[l] & & P_{i+1}=\vert L_{i+1} \vert \ar[ll]^{f_i^{i+1}} & \cdots\ar[l]   & X
}
\end{displaymath}

\vspace{1mm}

The existence of $\pi:Z \ra X$ with the required properties of fibers follows from Lemma~\ref{genWa}, when $P_i=\vert L_i\vert$,
$M_i=\vert L_i^{(n)}\vert$ and $\phi_i$ is the inclusion $i: \vert L_i^{(n)}\vert \hookrightarrow \vert L_i\vert$.

Note that $\phi_i^{-1}(B_{\da(i)}(x_i))=B_{\da(i)}(x_i)\cap \vert L_i^{(n)}\vert$, so (IV) of Lemma~\ref{genWa} becomes:

\begin{enumerate}
\item[(IV$^\ast$)]$\pi^{-1}(x) = \pi^{-1}((x_i))  =  \lim
(B_{\da(i)}(x_i)\cap \vert L_i^{(n)} \vert, g_i^{i+1})\ =\ \lim
(B_{\e(i)}(x_i)\cap \vert L_i^{(n)}\vert , g_i^{i+1}).$
\end{enumerate}

Property (IV) will guarantee that, for any $x\in X$,
$\pi^{-1}(x)\neq \emptyset$. This is true because, if we take any
$x=(x_i)\in X$, $x_i\in P_i=\vert L_i\vert$ implies that there is
a simplex $\sa \in L_i$ such that $x_i \in \sa$. Since $\mesh
L_i<\e(i)$, we get that $\diam \sa<\e(i)$, so $\sa \subset
B_{\e(i)}(x_i)$. Therefore $\sa^{(n)} \subset B_{\e(i)}(x_i)\cap
\vert L_i^{(n)}\vert$, so $$\emptyset \neq B_{\e(i)}(x_i)\cap \vert L_i^{(n)}\vert \subset B_{\da(i)}(x_i)\cap \vert L_i^{(n)}\vert=\phi_i^{-1}(B_{\da(i)}(x_i)) .$$ By (V) of Lemma~\ref{genWa}, $\pi : Z \ra X$ is surjective.

\vspace{2mm}

It remains to show that properties (V) and (VI) imply that $\pi$
is cell-like. Note that from (V) and (IV$^\ast$) we get that $\pi^{-1}(x)\
=\lim \ (P_{x_i,i}\cap \vert L_i^{(n)} \vert, g_i^{i+1})$, where $g_i^{i+1}$
stands for the appropriate restriction. It will be sufficient to
show that the maps $g_i^{i+1}: P_{x_{i+1},i+1}\cap \vert L_{i+1}^{(n)}\vert
\ra P_{x_{i},i}\cap \vert L_{i}^{(n)}\vert$ are null-homotopic.

First note that $P_{x_{i+1},i+1}$ being contractible implies that
the inclusion map $i: P_{x_{i+1},i+1}\cap \vert L_{i+1}^{(n)}\vert
\hookrightarrow P_{x_{i+1},i+1}$ is null-homotopic.  Since $\dim
P_{x_{i+1},i+1}\cap \vert L_{i+1}^{(n)}\vert \leq n$, $i$ is null-homotopic
as a map into $P_{x_{i+1},i+1}\cap \vert L_{i+1}^{(n+1)}\vert$, that is, this
homotopy happens within the $(n+1)$-skeleton of $L_{i+1}$. This is because $\dim\ ((P_{x_{i+1},i+1}\cap \vert L_{i+1}^{(n)}\vert )\times I )\leq n+1$, so if $H:(P_{x_{i+1},i+1}\cap \vert L_{i+1}^{(n)}\vert )\times I \ra P_{x_{i+1},i+1}$ is our homotopy, then, by Remark~\ref{unstable}, in each cell of $P_{x_{i+1},i+1}=\vert L_{i+1}\vert$ with dimension $\geq n+2$, the map $H$ will have unstable values.

 Using the last part of Remark~\ref{unstable}, as well as properties of deformation retracts, we can find a map $\widetilde H:(P_{x_{i+1},i+1}\cap \vert L_{i+1}^{(n)}\vert )\times I \ra P_{x_{i+1},i+1}$ such that $\widetilde H\vert_{H^{-1}(\vert L_{i+1}^{(n+1)}\vert)}=H\vert_{H^{-1}(\vert L_{i+1}^{(n+1)}\vert)}$, $\widetilde H ((P_{x_{i+1},i+1}\cap \vert L_{i+1}^{(n)}\vert) \times I)\subset \vert L_{i+1}^{(n+1)}\vert$, and so that $\widetilde H$ is a  homotopy between $i: P_{x_{i+1},i+1}\cap \vert L_{i+1}^{(n)}\vert
\hookrightarrow P_{x_{i+1},i+1}\cap \vert L_{i+1}^{(n+1)}\vert$ and a constant map.

Composing such a homotopy
with $g_{i}^{i+1}\vert_{\vert
L_{i+1}^{(n+1)}\vert}:\vert L_{i+1}^{(n+1)}\vert \ra \vert
L_{i}^{(n)}\vert$ yields the sought after null-homotopy for the
restriction $g_i^{i+1}\vert_{P_{x_{i+1},i+1}\cap
\vert L_{i+1}^{(n)}\vert}$.\hfill $\square$

\vspace{5mm}

We will now prepare for Lemma~\ref{V-lemma}, which will be useful in the proof of the new version of Edwards' Theorem. First note the following:

\begin{remark} \label{rem}
\emph{
Each $k$-dimensional simplex is homeomorphic to $I^k$, so it is an absolute extensor for normal spaces, hence also for CW-complexes. In particular, for a simplex $\sa$ we have $\vert \sa\vert \ \tau  \vert \sa\vert$.}
\end{remark}




\vspace{1mm}

\begin{lemma} \label{k-sim nbhd}
Let $\sa$ be a $k$-dimensional simplex. Then there exists an open neighborhood $N$ of $\vert\partial \sa\vert$ in $\vert \sa \vert$, and a surjective map $s:\vert\sa\vert\ra\vert\sa\vert$ such that $s(N)\subset\vert\partial \sa\vert$, and $s\vert_{\vert\partial \sa\vert}=id$.
\end{lemma}

\noindent\textit{Proof:}
It suffices to prove the lemma in the case when $\sa=\Delta\subset \R^{k+1}$ is the standard $k$-dimensional simplex.
Consider the homothety $h_{B,\frac{1}{2}}:\Delta \ra \Delta$, centered at the barycenter $B$ of $\Delta$ with scale $\frac{1}{2}$, that is, every point $P\in\Delta$ is mapped to $h_{B,\frac{1}{2}}(P)$ so that $B- h_{B,\frac{1}{2}}(P) =\frac{1}{2}(B-P)$.
Since $h_{B,\frac{1}{2}}(\Delta)$ is contained in the interior of $\Delta$, we see that $N:= \Delta \setminus h_{B,\frac{1}{2}}(\Delta)$ is an open neighborhood of $\partial\Delta$. Let $s:\Delta \ra\Delta$ be the map which on $h_{B,\frac{1}{2}}(\Delta)$ coincides with $(h_{B,\frac{1}{2}})^{-1}$, and on $N$ coincides with the restriction to $N$ of the central projection $\Delta\setminus B \ra \partial \Delta$.

\hfill $\square$

\vspace{3mm}

Using the previous Lemma we get the following
technical result helpful in the proof of Lemma~\ref{V-lemma}:

\begin{lemma} \label{tech}
Let $C$ be a finite simplicial complex with $\dim C =q$. Then for each $0\leq k \leq q$ there is an open neighborhood $U$ of
$\vert C^{(k)}\vert$ in $\vert C\vert$, and a surjective map $r:\vert
C\vert \ra \vert C\vert$ so that
\begin{enumerate}
\item[(1)] $r\vert_{\vert C^{(k)}\vert} = id_{\vert
C^{(k)}\vert}$,

\item[(2)] $r$ preserves simplexes, i.e., for any $\tau \in C$,
$r(\tau) \subset \tau$, and

\item[(3)] $r(U)\subset \vert C^{(k)}\vert$.
\end{enumerate}
\end{lemma}

\noindent \textit{Proof}:
The statement of this Lemma is true when $q=0$.
If $q\geq 1$ and $k=q-1$, the statement can be easily proven using Lemma~\ref{k-sim nbhd}.

Assume that $q>1$, and assume inductively that the statement of this Lemma is true when $q$ is replaced by $n$, and $0\leq n<q$.

Choose an open neighborhood $M$ of $\vert C^{(q-1)}\vert$ in $\vert C\vert$, and a surjective map $p:\vert C\vert \ra\vert C\vert$ so that
\begin{enumerate}
\item[(1)$_{q-1}$] $p\vert_{\vert C^{(q-1)}\vert} = id_{\vert C^{(q-1)}\vert}$,

\item[(2)$_{q-1}$]  $p(\tau) \subset \tau$ for any $\tau \in C$, and

\item[(3)$_{q-1}$] $p(M)\subset \vert C^{(q-1)}\vert$.
\end{enumerate}

If $k=q-1$, put $U:=M$ and $r:=p$ and we are done.
If $k<q-1$, proceed as follows. By the inductive assumption, we may select an open neighborhood $N$ of $\vert C^{(k)}\vert$ in $\vert C^{(q-1)}\vert$, and a surjective map $s:\vert C^{(q-1)}\vert\ra \vert C^{(q-1)}\vert$ such that
\begin{enumerate}
\item[(a)] $s\vert_{\vert C^{(k)}\vert} = id_{\vert C^{(k)}\vert}$,

\item[(b)]  $s(\tau) \subset \tau$ for any $\tau \in C^{(q-1)}$, and

\item[(c)] $s(N)\subset \vert C^{(k)}\vert$.
\end{enumerate}

For each $q$-simplex $\sa$ of $C$, $s(\vert\partial \sa\vert)=\vert\partial \sa\vert$ and $s\vert_{\vert\partial \sa\vert}:\vert\partial\sa\vert \ra \vert\partial\sa\vert$ is homotopic to identity. Hence there is a map $s_\sa:\vert\sa\vert\ra\vert\sa\vert$ such that $s_\sa\vert_{\vert\partial\sa\vert}=s\vert_{\vert\partial\sa\vert}$, and $s_\sa$ must be surjective.
Put $\widetilde s:=s\ \cup \ \left(\cup \ \{s_\sa \vert \ \sa \ \mathrm{is \ a \ } q\mathrm{-simplex \ of \ }C\}\right)$.
Then $\widetilde s:\vert C\vert \ra \vert C\vert $ is surjective, $\widetilde s(\tau) \subset \tau$ for any $\tau \in C$, and
$\widetilde s\vert_{\vert C^{(q-1)}\vert}=s$.

Note that $p\vert_M:M\ra \vert C^{(q-1)}\vert$ is continuous, and $N$ is open in $\vert C^{(q-1)}\vert$, so $(p\vert_M)^{-1}(N)$ is open in $M$ and therefore open in $\vert C\vert$.

Define $U:=(p\vert_M)^{-1}(N)=M\cap p^{-1}(N)$ and $r:=\widetilde s\circ p:\vert C\vert \ra \vert C\vert$. Observe that $U$ is a neighborhood of $\vert C^{(k)}\vert$ in $\vert C\vert$ and that $r$ is surjective. It is routine to check that (1)--(3) are true.
\hfill $\square$

\vspace{3mm}


\begin{lemma}\label{V-lemma}
For any finite simplicial complex $C$, there is a map $r:\vert C\vert \ra \vert C\vert$ and an open cover
$\mathcal{V}=\{V_\sa : \sa \in C\}$ of $\vert C\vert$ such that for all $\sa$, $\tau \in C$:
\begin{enumerate}
\item[(i)] $\overset{\circ}\sa \subset V_\sa$,

\item[(ii)] if $\sa\neq \tau$ and $\dim \sa = \dim \tau$, $V_\sa$ and $V_\tau$ are disjoint,

\item[(iii)] if $y\in \overset{\circ}\tau$,
 $\dim \sa \geq \dim \tau$ and $\sa \neq \tau$, then  $y \notin V_\sa$,

\item[(iv)]if $y \in \overset{\circ}\tau \cap V_\sa$, where  $\dim \sa <\dim \tau$, then $\sa$ is a face of $\tau$, and

\item[(v)] $r(V_\sa)\subset \sa$.
\end{enumerate}
\end{lemma}

\noindent \textit{Proof}:
Since $C$ is finite, let us suppose that $\dim C = q$.
%
%
For $k=0,\ldots ,q-1$, let $U_k$ correspond to $U$ and $r_k$ correspond to $r$ from Lemma~\ref{tech}.
Note that for vertices $v
\in C^{(0)}$ we have that $\overset{\circ}v = v$.

\vspace{2mm}

Here is how we will define the open cover $\mathcal{V}=\{V_\sa : \sa \in C\}$ for $\vert
C\vert$:
\begin{enumerate}
\item[(a)] for each $k$-simplex $\sa$ of $C$, where $k=0,\ldots
,q-1$, put\\ $V_{\sa}:=\ (r_k\circ r_{k+1}\circ \ldots \circ
r_{q-1})^{-1}(\overset{\circ}\sa)$ into $\mathcal{V}$, and

\item[(b)] for each $q$-simplex $\sa$ of $C$, put
$V_{\sa}:=\overset{\circ}\sa$ into $\mathcal{V}$.
\end{enumerate}

Note that all elements of $\mathcal{V}$ are open sets: in (b) that
is clear, and in (a): $(r_k\circ r_{k+1}\circ \ldots \circ
r_{q-1})^{-1}(\overset{\circ}\sa)=\ r_{q-1}^{-1}(\ldots
(r_{k+1}^{-1}(r_k^{-1}(\overset{\circ}\sa))))$, and
$r_k^{-1}(\overset{\circ}\sa)$ is open because $r_k\vert_{U_k}:U_k
\ra \vert C^{(k)}\vert$ is continuous, and $\overset{\circ}\sa$ is open in
$\vert C^{(k)}\vert$.

Let us check that (i) is true: $\overset{\circ}\sa \subset V_\sa$
is clear for case (b), and, for case (a), since $r_k,
r_{k+1},\ldots , r_{q-1}$ are all the identity on $\vert
C^{(k)}\vert$ and $\overset{\circ}\sa \subset \vert
C^{(k)}\vert$, then $\overset{\circ}\sa \subset V_\sa$. Hence
$\mathcal{V}$ is a cover for $\vert C\vert$ because of (i).

If $\sa$ and $\tau$ are two different
simplexes of the same dimension, then $\overset{\circ}\sa$ and
$\overset{\circ}\tau$ are disjoint. If $\dim \sa = \dim \tau = q$, (ii) is clear.
If $\dim \sa = \dim \tau < q$, then (a) implies that $V_{\sa}$
and $V_{\tau}$ are disjoint, i.e., (ii) is true.

Let us prove property (iii). We know that $y\in\overset{\circ}\tau
\subset V_\tau$. If $\tau$ and $\sa$ are of the same dimension,
then (ii) implies $y\notin V_\sa$. If $\dim
\tau <\dim \sa \leq q-1$, then $V_{\sa}:=\ (r_{\dim \sa}\circ
\ldots \circ r_{q-1})^{-1}(\overset{\circ}\sa)$, so if $y$ would
be in $ V_\sa$, then $r_{\dim \sa}\circ \ldots \circ r_{q-1}(y)\in
\overset{\circ}\sa$. But $r_{\dim \sa}, \ldots , r_{q-1}$ are the
identity on $\vert C^{(\dim \tau)}\vert \supset \tau$, so $r_{\dim \sa}\circ
\ldots \circ r_{q-1}(y)=y \in \overset{\circ}\sa$, which is in
contradiction with $y \in \overset{\circ}\tau$. Thus $y\notin
V_\sa$. If $\dim \tau <\dim \sa = q$, then
$V_\sa=\overset{\circ}\sa$, so $y\in \overset{\circ}\tau$ and
$\tau \neq \sa$ imply that $y\notin V_\sa$.

To prove (iv), suppose that $y \in V_\sa$
for some $\sa \in C$ with $\dim \sa < \dim \tau$. Then $V_{\sa}:=\
(r_{\dim \sa}\circ \ldots \circ
r_{q-1})^{-1}(\overset{\circ}\sa)$, so $r_{\dim \sa}\circ \ldots
\circ r_{q-1}(y)\in \overset{\circ}\sa$. Notice that $r_{\dim
\tau},r_{\dim \tau +1}, \ldots , r_{q-1}$ are the identity on $\tau$,
so $r_{\dim \sa}\circ \ldots \circ r_{q-1}(y)=r_{\dim \sa}\circ
\ldots \circ r_{\dim \tau -1}(y) \in \overset{\circ}\sa$. The maps
$r_{\dim \sa}, \ldots , r_{\dim \tau -1}$ preserve simplexes, by
(2) of Lemma~\ref{tech}, so $y \in \overset{\circ}\tau$ implies that $r_{\dim
\sa}\circ \ldots \circ r_{\dim \tau -1}(y) \in \tau$. Thus $\tau
\cap \overset{\circ}\sa \neq \emptyset$, so $\sa$ must be a face
of $\tau$.

It remains to define the map $r$ and prove the property (v). Define $r:=r_0\circ r_1\circ\ldots\circ r_{q-1}:\vert C\vert \ra \vert C\vert$. For any $k$-simplex $\sa$ of $C$ where $k=1,\ldots
,q-1$, by (a) we get that
$$
r(V_\sa) = r_0\circ r_1\circ\ldots \circ
r_{q-1}((r_k\circ r_{k+1}\circ \ldots \circ r_{q-1})^{-1}(\overset{\circ}\sa))\\
= r_0\circ r_{1}\circ \ldots \circ
r_{k-1}(\overset{\circ}\sa),
$$
 since all $r_i$ are surjective. Also, by (2) of Lemma~\ref{tech}, $r(V_\sa)=r_0\circ r_{1}\circ \ldots \circ r_{k-1}(\overset{\circ}\sa)\subset \sa$.

 Likewise, for any $q$-simplex $\sa$ of $C$, we get
$r(V_\sa)=r(\overset{\circ}\sa)\subset \sa$ for the same reason.
For vertices $v\in C^{(0)}$, $r(V_v)= r\circ r^{-1}(v)=v$.
So we conclude that (v) is true. \hfill $\square$

\vspace{5mm}

\noindent Next we will see a version of Theorem 4.2 from \cite{Wa}, adapted for our
situation. In order to proceed, however, we will need to be reminded of two definitions.

Let $K$ be a simplicial complex,  $X$  a space, and  $f:X \ra \vert K\vert$  a map. Recall that a map $g:X \ra \vert K\vert$ is called a $K$-\emph{modification} of $f$ if whenever $x\in X$ and $f(x)\in \sa$, for some $\sa \in K$, then $g(x)\in \sa$. This is equivalent to the following: whenever $x\in X$ and $f(x)\in \overset{\circ}\sa$, for some $\sa \in K$, then $g(x)\in \sa$.

In the course of the proof of the following theorem, we will need the notion of \textit{resolution in the sense of inverse sequences}. This usage of the word resolution is completely different from the notion from the title of this paper. The definition can be found in \cite{MS1} for the more general case of inverse systems. Here, however, we will give the definition for inverse sequences.

\begin{defi} \label{invsecrez}
\emph{
Let $X$ be a topological space. A \emph{resolution} of $X$ \emph{in the sense of inverse sequences} consists of an inverse sequence of topological spaces $\mathbf{X}= (X_i,p_i^{i+1})$ and a family of maps $(p_i:X \ra X_i)$ with the following two properties:
\begin{enumerate}
\item[(R1)] Let $P$ be an ANR, $\scv$ an open cover of $P$ and $h:X\ra P$ a map. Then there is an index $s\in \N$ and a map $f:X_s \ra P$ such that the maps $f\circ p_s$ and $h$ are $\scv$-close.
\item[(R2)] Let $P$ be an ANR and $\scv$ an open cover of $P$. There exists an open cover $\scv '$ of $P$ with the following property: if $s\in\N$ and $f,f':X_s \ra P$ are maps such that the maps $f\circ p_s$ and $f'\circ p_s$ are $\scv '$-close, then there exists an $s'\geq s$ such that the maps $f\circ p_s^{s'}$ and $f'\circ p_s^{s'}$ are $\scv$-close.
\end{enumerate}}
\end{defi}

By Theorem I.6.1.1 from \cite{MS1}, if all $X_i$ in $\mathbf{X}$ are compact Hausdorff spaces, then $\mathbf{X}= (X_i,p_i^{i+1})$ with its usual projection maps $(p_i:\lim X \ra X_i)$ is a resolution of $\lim X$ in the sense of inverse sequences.

Moreover, since every compact metrizable space $X$ is the inverse limit of an inverse sequence of compact polyhedra $\mathbf{X}= (P_i,p_i^{i+1})$ (see Corollary I.5.2.4 of \cite{MS1}), this inverse sequence $\mathbf{X}$ will have the property (R1) mentioned above, and we will refer to this property as the \emph{resolution property} (R1) \emph{in the sense of inverse sequences}.

\begin{tm}[New statement of Edwards Theorem] \label{Ed} Let $n\in \N$ and let $Y$ be a
compact metrizable space such that $Y=\lim \ (\vert L_i\vert,
f_i^{i+1})$, where $\vert L_i\vert$ are compact polyhedra with
$\dim L_i \leq n+1$, and $f_i^{i+1}$ are surjections. Then
$\dim_\Z Y \leq n$ implies that there exists an $s\in \N$, $s >1$,
and there exists a map $g_1^s : \vert L_s\vert \to \vert
L_1^{(n)}\vert$ which is an $L_1$-modification of $f_1^s$.
\end{tm}
\begin{displaymath}
\xymatrix{
\vert L_1^{(n)}\vert \ar@{_{(}->}[d] & &                                  &               &\\
 \vert L_1\vert &  & \vert L_s\vert \ar[ll]^{f_1^s}\ar@{-->}[llu]_{g_1^s} & \cdots \ar[l] & Y
}
\end{displaymath}

\vspace{2mm}

\noindent \textit{Proof}: There will be two separate parts of this
proof, for $n\geq 2$ and for $n=1$.

\vspace{2mm}

Let us start with $n\geq 2$. We will build an Edwards-Walsh
complex $\widehat{L}_1$ above $L_1^{(n)}$. Since $\dim L_1\leq
n+1$ and $L_1$ is finite, $L_1$ has to have finitely many
$(n+1)$-simplexes, say, $\sa_1,\ldots , \sa_m$. Focus on
$L_1^{(n)}$, and above each of $\sa_i^{(n)}=\partial
\sa_i\thickapprox S^n$, build a $K(\Z,n)$ by attaching cells of
dimension $(n+2)$ and higher. Name the CW-complex that we get in
this fashion $\widehat{L}_1$. Notice that we can write
$\widehat{L}_1=L_1^{(n)} \cup K(\sa_1)\cup K(\sa_2)\cup \ldots
\cup K(\sa_m)$, where each $K(\sa_i)$ is a $K(\Z,n)$ attached to
$\partial \sa_i$. Also notice that we can make the attaching maps
piecewise linear, so that we will be able to triangulate
$\widehat{L}_1$ keeping $L_1^{(n)}$ as a subcomplex.

Let $\theta:\widehat{L}_1 \ra \vert L_1\vert$ be a map such that
$\theta\vert_{\vert L_1^{(n)}\vert}=id_{\vert L_1^{(n)}\vert}$ and
$\theta (K(\sa_i))\subset \sa_i$. This $\theta$ can be constructed
as follows: first, define $\theta\vert_{\vert
L_1^{(n)}\vert}:=id_{\vert L_1^{(n)}\vert}$. By Remark~\ref{rem}, each $\sa_i$ is an absolute extensor for CW-complexes,
so the inclusion map $j:\sa_i^{(n)} \ra \sa_i$ can be extended over
$K(\sa_i)$. Call this extension $\theta\vert_{K(\sa_i)}$. Gluing
together all of the extensions $\theta\vert_{K(\sa_i)}$ for
$i=1,\ldots ,m$ with $\theta\vert_{\vert L_1^{(n)}\vert}$ will
produce the map $\theta$.

Let $f_1:Y \ra \vert L_1\vert$ be the projection map from the
inverse sequence. The map $f_1$ is surjective since all  $f_i^{i+1}$
are surjective. Extend $f_1\vert_{f_1^{-1}(\vert L_1^{(n)}\vert)}:
f_1^{-1}(\vert L_1^{(n)}\vert)\ra \vert L_1^{(n)}\vert$ to a map
$h:Y\ra \widehat{L}_1$ such that

\begin{enumerate}
\item[(a)] $h(f_1^{-1}(\sa_i))\subset
\theta^{-1}(\sa_i)=K(\sa_i)$, for $i=1,\ldots ,m$.
\end{enumerate}

This can be
done using $\dim_\Z Y\leq n \ \Leftrightarrow \ Y\tau K(\Z,n)$:
for any $(n+1)$-dimensional $\sa_i$, take
$f_1\vert_{f_1^{-1}(\sa_i^{(n)})}:f_1^{-1}(\sa_i^{(n)})\ra
\sa_i^{(n)}$ and compose it with the inclusion
$i:\sa_i^{(n)}\hookrightarrow K(\sa_i)=K(\Z,n)$. Now $Y\tau K(\Z,n)$
implies $f_1^{-1}(\sa_i) \ \tau K(\Z,n)$, so the map $i\circ
f_1\vert_{f_1^{-1}(\sa_i^{(n)})}: f_1^{-1}(\sa_i^{(n)})\ra
K(\sa_i)$ can be extended over $f_1^{-1}(\sa_i)$. Call this
extension $h\vert_{f_1^{-1}(\sa_i)}$. So we get the map $h$ that
we need by gluing together all of the extensions
$h\vert_{f_1^{-1}(\sa_i)}$, for $i=1,\ldots ,m$, with
$h\vert_{f_1^{-1}(\vert L_1^{(n)}\vert)}=f_1\vert_{f_1^{-1}(\vert
L_1^{(n)}\vert)}$.

Note that our inverse sequence $(\vert L_i\vert, f_i^{i+1})$ is a
compact resolution for $Y$ in the sense of inverse sequences (see Definition~\ref{invsecrez}), so, in particular, it has the
resolution property (R1) (in the sense of inverse sequences): if we choose an open cover $\mathcal{V}$
for the minimal and hence finite subcomplex  $\widehat{C}$ in $\widehat{L}_1$ such
that $h(Y)\subset \widehat{C}$, then we can find an $s>1$ and a
map $h_1^s:\vert L_s\vert \ra \widehat{C}$ such that $h$ and
$h_1^s\circ f_s$ are $\scv$-close.

\begin{displaymath}
\xymatrix{
 \widehat{L}_1  \ar[dd]_{\theta}& \ \widehat{C}\ar@{_{(}->}[l] & & & &\\
                                                               & & & & \\
\vert L_1\vert & & \vert L_s\vert \ar[ll]^{f_1^s} \ar@{-->}[uul]_{\hspace{-2mm} h_1^s} \ar@/^/@{.>}[uul]|{h_s}& ...\ar[l] & & Y \ar@/^/[lll]^{\qquad f_s}  \ar[uullll]_h
          }
\end{displaymath}

\vspace{2mm}
Let us make a wise choice for $\mathcal{V}$. Start by
triangulating $\widehat{C}$: let $C$ denote a finite simplicial
complex which is a triangulation of $\widehat{C}$ whose restriction to
$\vert L_1^{(n)}\vert$ is a subcomplex. So $\vert C\vert = \widehat{C}$.
Since $C$ is finite, let us suppose that $\dim C = q$.

Define an open cover $\mathcal{V}$ for $\vert C \vert$, and a map
$r: \vert C\vert \ra \vert C\vert$ as in Lemma~\ref{V-lemma}. For
this cover $\mathcal{V}$ for $\vert C\vert$, we may apply
resolution property (R1) (in the sense of inverse sequences):  we can find an $s>1$ and a map
$h_1^s:\vert L_s\vert \ra \vert C\vert$ such that $h$ and
$h_1^s\circ f_s$ are $\mathcal{V}$-close. Define $h_s:=r\circ
h_1^s:\vert L_s\vert \ra \vert C\vert.$ Because of our choices, we
get that
\begin{enumerate}
\item[(b)] whenever $h(y)\in \overset{\circ}\tau$ for some $\tau
\in C$, then $(h_s\circ f_s) (y) \in \tau$.
\end{enumerate}
This is true because, by (i), (ii), (iii) and (iv) of Lemma~\ref{V-lemma}, $h(y)\in \overset{\circ}\tau$
implies that $h(y) \in V_\tau$, and possibly also $h(y) \in V_\sa$ for
some $\sa$ which is a face of $\tau$, but $h(y)$ is in no other elements of $\mathcal{V}$. Since $h_1^s\circ
f_s$ is $\mathcal{V}$-close to $h$, we have that either
$h_1^s\circ f_s(y)\in V_\tau$, or $h_1^s\circ f_s(y)\in V_\sa$,
for some face $\sa$ of $\tau$. But by (v) of Lemma~\ref{V-lemma}, $r(V_\tau) \subset
\tau$ and $r(V_\sa) \subset \sa \subset \tau$. Thus $h_s\circ
f_s(y)=r\circ h_1^s\circ f_s(y)\in \tau$.

If $f_1(y)\in \sa_i$ for some $(n+1)$-simplex $\sa_i$ of $L_1$,
then, by (a), $h(y) \in K(\sa_i)$, so $h(y)\in \overset{\circ}\tau$ for
some $\tau \in C$ and $\tau\subset K(\sa_i)$. By (b), $h_s(f_s(y))\in \tau$. So
we can conclude that
\begin{enumerate}
\item[(c)] if $f_1(y)\in \sa_i$, for some $(n+1)$-simplex $\sa_i$
of $L_1$, then both $h(y)$ and $h_s\circ f_s(y)$ land in
$K(\sa_i)$.
\end{enumerate}

Now we will construct a map $g_1^s:\vert L_s\vert \ra \vert
L_1^{(n)}\vert$ such that :
\begin{enumerate}
\item[(d)] $g_1^s\vert_{h_s^{-1}(\vert
L_1^{(n)}\vert)}=h_s\vert_{h_s^{-1}(\vert L_1^{(n)}\vert)}$, and

\item[(e)] whenever $h_s(z) \in K(\sa_i)$ for some
$(n+1)$-simplex $\sa_i$ of $L_1$, then $g_1^s(z)\in \sa_i$.
\end{enumerate}
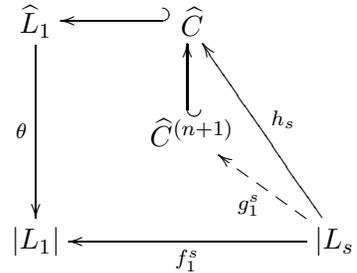
\begin{displaymath}
\xymatrix{
 \widehat{L}_1  \ar[dd]_{\theta}& \ \widehat{C}\ar@{_{(}->}[l] &  \\
                                &  \ \widehat{C}^{(n+1)} \ar@{_{(}->}[u]    &  \\
\vert L_1\vert & & \vert L_s\vert \ar[ll]^{f_1^s} \ar[uul]_{h_s} \ar@{-->}[ul]^{g_1^s}
          }
\end{displaymath}
We know that $h_s: \vert L_s\vert \ra \vert C\vert=\widehat{C}$,
where $C$ is a triangulation of the finite CW-subcomplex
$\widehat{C}$ of $\widehat{L}_1$. Since $\widehat{C}$
is finite, we can pick a cell $\ga$ of maximal possible dimension
$\dim \ga = q$ (we have assumed that $\dim C =q$, so $\dim
\widehat{C}=q$). It is safe to assume that $q\geq n+2$.

Pick a point $w$ in $\overset{\circ}\ga$ with an open neighborhood
$W \subset \overset{\circ}\ga$. Since $\dim \vert L_s\vert \leq
n+1$ and $\dim \ga >n+1$, the point $w$ we picked is an unstable
value for $h_s$, so we can construct a new map
$g_{1,\ga}^s: \vert L_s \vert \ra \widehat{C}\setminus \{w\}$
that agrees with $h_s$ on $h_s^{-1}(\widehat{C}\setminus W)$,
and $g_{1,\ga}^s(h_s^{-1}(\ga))\subset \ga \setminus \{w\}$  (see Remark~\ref{unstable}).
Retract $\ga \setminus \{w\}$ to $\partial \ga$ by a retraction
$\tilde{r}:\widehat{C}\setminus \{w\} \ra \widehat{C}\setminus
\overset{\circ}\ga$, such that
$\tilde{r}\vert_{\widehat{C}\setminus \overset{\circ}\ga}=id$.
Replace $h_s$ with $\tilde{r}\circ g_{1,\ga}^s: \vert L_s\vert \ra
\widehat{C}\setminus \overset{\circ}\ga$.

We will repeat this process, starting with $\widehat{C}\setminus
\overset{\circ}\ga$ and the map $\tilde{r}\circ g_{1,\ga}^s$
instead of $\widehat{C}$ and $h_s$: pick a cell of maximal
dimension in $\widehat{C}\setminus \overset{\circ}\ga$, etc. This
is done one cell at a time, until we get rid of all cells in
$\widehat{C}$ with dimension $\geq n+2$. The map we end up with
will be $g_1^s:\vert L_s\vert \ra \widehat{C}^{(n+1)}$, where
$\widehat{C}^{(n+1)}$ stands for the CW-skeleton of dimension $n+1$
for $\widehat{C}$. Notice that $\widehat{C}^{(n+1)}\subset
\widehat{L}_1^{(n+1)}$, but the CW-skeleton of dimension $n+1$ for
$\widehat{L}_1$ is equal to the CW-skeleton of dimension $n$ for
$\widehat{L}_1$, since we have built $\widehat{L}_1$ by attaching
cells of dimension $n+2$ and higher to $L_1^{(n)}$. Thus
$\widehat{L}_1^{(n+1)}=\widehat{L}_1^{(n)}=\vert L_1^{(n)}\vert$,
where $L_1^{(n)}$ is the simplicial $n$-skeleton of $L_1$. So in
fact, $g_1^s:\vert L_s\vert \ra \vert L_1^{(n)}\vert$.

By our construction, $g_1^s$ agrees with $h_s$ on $h_s^{-1}(\vert
L_1^{(n)}\vert)$, so (d) is true. To prove property (e), let
$h_s(z) \in K(\sa_i)$. Then $h_s(z) \in \ga$, for some cell $\ga$
of $K(\sa_i)$. So $\tilde{r}\circ g_{1,\ga}^s(z) \in \partial \ga
\subset K(\sa_i)$. As we go on with our construction, we get
$g_1^s(z)\in (K(\sa_i))^{(n+1)}=\partial \sa_i\subset \sa_i$.

\vspace{2mm}

Finally, for any $z \in \vert L_s\vert$ we have that either
$f_1^s(z) \in \overset{\circ}\tau$, for some $\tau \in L_1^{(n)}$,
or $f_1^s(z) \in \overset{\circ}\sa_i$, for some $(n+1)$-simplex
$\sa_i$ of $L_1$. Since $f_s$ is surjective, there is a $y \in Y$
such that $f_s(y)=z$.

So, if $f_1^s(z) \in \overset{\circ}\tau$ for some $\tau \in
L_1^{(n)}$, then $f_1(y)=f_1^s(f_s(y))=f_1^s(z) \in
\overset{\circ}\tau \subset \vert L_1^{(n)}\vert$. Recall that on
$f_1^{-1}(\vert L_1^{(n)}\vert)$, $f_1$ and $h$ coincide. Thus
$f_1(y)=h(y) \in \overset{\circ}\tau$. There is a simplex
$\tau'\in C\cap \vert L_1^{(n)}\vert$ such that $\tau'\subset
\tau$, and $f_1(y)=h(y) \in {\overset{\circ}\tau} \ \!'$. By (b) we
get that $h_s\circ f_s(y) \in \tau' \subset \tau$, i.e., $h_s(z)
\in \tau \in L_1^{(n)}$, so by (d), $g_1^s(z)=h_s(z) \in \tau$.

On the other hand, if $f_1^s(z) \in \overset{\circ}\sa_i$, for
some $(n+1)$-simplex $\sa_i$ of $L_1$, then
$f_1(y)=f_1^s(f_s(y))=f_1^s(z) \in \overset{\circ}\sa_i$. By (c),
$h_s\circ f_s(y) \in K(\sa_i)$, i.e., $h_s(z) \in K(\sa_i)$.
Property (e) implies that $g_1^s(z) \in \sa_i$.

So $g_1^s$ is an $L_1$-modification of $f_1^s$.

\vspace{5mm}

It remains to prove this theorem for $n=1$. First note that
$\dim_\Z Y \leq 1$ implies that $\dim Y\leq 1$, because $S^1$ is a $K(\Z,1)$-complex. We will not need
to construct an Edwards-Walsh complex $\widehat{L}_1$ here.
Instead, look at the map $f_1:Y \ra \vert L_1\vert$. Let $g_1:Y
\ra \vert L_1^{(1)}\vert$ be a stability theory version of $f_1$.
We construct $g_1$ as before: since we know that $\dim L_1 \leq
2$, pick any $2$-simplex $\sa$ of $L_1$. We can pick a point $w\in
\overset{\circ}\sa$ with an open neighborhood $W \subset
\overset{\circ}\sa$, and since $\dim \sa =2$, the point $w$ is an
unstable value for $f_1$. So there exists a map $g_{1,\sa}:Y\ra \vert
L_1\vert \setminus \{w\}$ which agrees with $f_1$ on
$f_1^{-1}(\vert L_1\vert \setminus W)$, and such that
$g_{1,\sa}(f_1^{-1}(\sa))\subset \sa \setminus \{w\}$. Now retract
$\sa \setminus \{w\}$ to $\partial \sa$ by a retraction
$\tilde{r}$ which is the identity on $\vert L_1 \vert \setminus
\overset{\circ}\sa$. Finally, replace $f_1$ by $\tilde{r}\circ
g_{1,\sa}: Y \ra \vert L_1\vert \setminus \overset{\circ}\sa$.
Continue the process with one $2$-simplex at a time. Since $L_1$
is finite, in finitely many steps we will reach the needed map
$g_1:Y \ra \vert L_1^{(1)}\vert$. Note that from the construction of
$g_1$, we get
\begin{enumerate}
\item[(f)] $g_1\vert_{f_1^{-1}(\vert
L_1^{(1)}\vert)}=f_1\vert_{f_1^{-1}(\vert L_1^{(1)}\vert)}$, and
for every $2$-simplex $\sa$ of $L_1$, $\ g_1(f_1^{-1}(\sa))\subset
\partial \sa$.
\end{enumerate}
\begin{displaymath}
\xymatrix{
 \vert L^{(1)}_1 \vert \ar@{_{(}->}[d] & & & &\\
\vert L_1\vert &  & \vert L_s\vert \ar[ll]^{f_1^s} \ar@{-->}[ull]_{\hspace{-3mm}\widehat{g}_1^s} \ar@/^/@{.>}[ull]|{g_1^s} & ...\ar[l]  & Y \ar@/^/[ll]^{\qquad f_s}  \ar[ullll]_{g_1}
}
\end{displaymath}
Let us choose an open cover $\mathcal{V}$ of $L_1^{(1)}$ as
before: apply Lemma~\ref{V-lemma} to $C=L_1^{(1)}$. Note that $q=1$, so the map $r = r_0:\vert L_1^{(1)}\vert \ra \vert L_1^{(1)}\vert$.

Now we can use  resolution property (R1) (in the sense of inverse sequences): there is an index $s>1$
and a map $\widehat{g}_1^s:\vert L_s \vert \ra \vert
L_1^{(1)}\vert$ such that $\widehat{g}_1^s\circ f_s$ and $g_1$ are
$\mathcal{V}$-close. Define $g_1^s:=r_0\circ \widehat{g}_1^s
:\vert L_s \vert \ra \vert L_1^{(1)}\vert$.

Notice that for any $y\in Y$, if $g_1(y)\in \overset{\circ}\tau$
for some $\tau \in L_1^{(1)}$ (vertices included), then $g_1(y)\in
V_\tau$, and possibly also $g_1(y) \in V_v$, where $v$ is a vertex
of $\tau$. Then either $\widehat{g}_1^s \circ f_s(y) \in V_\tau$,
or $\widehat{g}_1^s \circ f_s(y) \in V_v$. In any case, $r_0\circ
\widehat{g}_1^s \circ f_s(y) \in \tau$. Hence,
\begin{enumerate}
\item[(g)] for any $y\in Y$, $g_1(y)\in \overset{\circ}\tau$ for
some $\tau \in L_1^{(1)}$, implies that $g_1^s(f_s(y)) \in \tau$.
\end{enumerate}

Finally, for any $z \in \vert L_s\vert$, $f_s$ is surjective
implies that there is a $y \in Y$ such that $f_s(y)=z$. Then
$f_1^s(z)=f_1^s(f_s(y))=f_1(y)$. Now $f_1^s(z)$ is either in
$\overset{\circ}\sa$ for some $2$-simplex $\sa$ in $L_1$, or in
$\overset{\circ}\tau$ for some $\tau \in L_1^{(1)}$.

If $f_1^s(z)\in \overset{\circ}\sa$, that is $f_1(y) \in
\overset{\circ}\sa$ for some $2$-simplex $\sa$, by (f) we get
that $g_1(y)\in
\partial \sa$. Then by (g), $g_1^s(f_s(y)) \in \partial \sa$,
i.e., $g_1^s(z) \in \sa$.

If $f_1^s(z)=f_1(y)\in \overset{\circ}\tau$ for some $\tau \in
L_1^{(1)}$, then (f) implies that $g_1(y)=f_1(y) \in
\overset{\circ}\tau$, so by (g), $g_1^s(f_s(y)) \in \tau$, i.e.,
$g_1^s(z)\in \tau$.

Therefore, $g_1^s$ is indeed an $L_1$-modification of
$f_1^s$.\hfill $\square$

\vspace{2mm}

\begin{lemma}
\label{dimGY} Let $n\in \N$, $G$ be an abelian group and
$K$ be a connected \emph{CW}-complex with $\pi_n(K)\cong G$,
$\pi_k(K)\cong 0$ for $0\leq k< n$. If $Y$ is a compact metrizable
space with $\dim Y\leq n+1$, then  $Y\tau K \ \Leftrightarrow\
\dim_G Y \leq n$.
\end{lemma}


\noindent \textit{Proof}: Build a $K(G,n)$ by attaching cells of
dimension $n+2$ and higher to our CW-complex $K$.

First assume that $Y\tau K$, and let us show that $\dim_G Y \leq n$.
If we look at any closed set $A\subset Y$ and any map $f:A \ra
K(G,n)$, we have that $\dim A \leq \dim Y \leq n+1$, so we can
homotope $f$ into $K(G,n)^{(n+1)}=K^{(n+1)}\subset K$, i.e., there
is a map $\overline{f} :A \ra K$ which is homotopic to $f$. Now $Y\tau K$ implies the existence of a map $g: Y \ra K$ which
 extends $\overline{f}$. Therefore, by the homotopy extension theorem,
 $f$ can be extended continuously over $Y$, so we get that $Y\tau K \ \Rightarrow \ Y \tau K(G,n) \ \Rightarrow \ \dim_G Y \leq
 n$.

Second, assume that $\dim_G Y \leq n$, and let us show $Y \tau K$.
Look at any closed set $A\subset Y$ and any map $f:A \ra K$. Let
$i:K \hookrightarrow K(G,n)$
be the inclusion map. Then $Y\tau K(G,n)$ implies that there is a map $\tilde{f}:Y \ra K(G,n)$
extending $i\circ f :A \ra K(G,n)$, i.e., $\tilde{f}\vert_A
= i \circ f$.

Since $Y$ is compact, $\tilde{f}(Y)$ is contained in a finite subcomplex $\widehat{C}$ of $K(G,n)$.
There are finitely many cells in $\widehat{C}\setminus K$, and all of them have dimension $\geq n+2$.
Pick a cell of maximal dimension $\ga \in \widehat{C}\setminus K$, and a point $w\in \overset{\circ}\ga$ with an open neighborhood $W\subset \overset{\circ}\ga$. Since $\dim Y \leq n+1$ and $\dim \ga \geq n+2$, the point $w$ is an unstable value of the map $\tilde{f}$,
so there is a map $g_\ga:Y \ra \widehat{C}\setminus \{w\}$ which agrees with $\tilde{f}$ on $\tilde{f}^{-1}(\widehat{C}\setminus W)$, and such that
$g_\ga(\tilde{f}^{-1}(\ga))\subset \ga \setminus \{w\}$. Retract $\ga \setminus \{w\}$ to $\partial \ga$ by a retraction
$\tilde{r}:\widehat{C}\setminus \{w\} \ra \widehat{C}\setminus
\overset{\circ}\ga$, such that
$\tilde{r}\vert_{\widehat{C}\setminus \overset{\circ}\ga}=id$.
Replace $\tilde{f}$ with $\tilde{r}\circ g_{\ga}: Y \ra
\widehat{C}\setminus \overset{\circ}\ga$. Repeat this process one cell at a time until all cells of $\widehat{C} \setminus K$ are exhausted. The map we end up with will be $g:Y \ra K$ such that $g\vert_{\tilde{f}^{-1}(K)}=\tilde{f}\vert_{\tilde{f}^{-1}(K)}$.
Since $\tilde{f}(A)=f(A)\subset K$, that is, $A\subset \tilde{f}^{-1}(K)$, we get $g\vert_A=\tilde{f}\vert_A$. So $g:Y \ra K$ is an extension of $f:A\ra K$. Therefore $Y\tau K$. \hfill $\square$

\section{Lemmas for inverse sequences}

The proof of the main result will require certain manipulations of inverse sequences of metric compacta. This section will contain the needed results, mostly taken from Section 3 of \cite{RS2}.
The next lemma follows from Corollary 1 of \cite{MS2}, or from \cite{Br}.

\begin{lemma}\label{stability}
Let $\mathbf{X}=(X_i,p_i^{i+1})$ be an inverse sequence of metric compacta $(X_i,d_i)$. Then there exists a sequence $(\ga_i)$ of positive numbers such that if $\mathbf{Y}=(X_i,q_i^{i+1})$ is an inverse sequence and $d_i(p_i^{i+1},q_i^{i+1})<\ga_i$ for each $i$, then $\lim \mathbf{Y}=\lim \mathbf{X}$.\hfill $\square$
\end{lemma}

We shall call such $(\ga_i)$ a \emph{sequence of stability} for $\mathbf{X}$.

\vspace{2mm}

Let $K$ be a simplicial complex,  $X$  a space, and  $f:X \ra \vert K\vert$  a map.
One calls $f$ a $K$-\emph{irreducible} map if each $K$-modification $g$ of $f$ is surjective. Note that $f$ being
$K$-irreducible implies that $f$ is surjective, and for any subdivision $M$ of $K$, $f$ is $M$-irreducible.

\begin{lemma}
\label{K-irred}
If $f:X \ra \vert K\vert$ is a $K$-irreducible map, and $g:X \ra \vert K \vert$ is a \linebreak $K$-modification of $f$, then $g$ is $K$-irreducible. \hfill $\square$
\end{lemma}

\vspace{3mm}

The following fact may be deduced from Theorem 3.11 of \cite{JR}, or found in \cite{Fr} (Hauptsatz I, p.\ 191).

\begin{lemma}
\label{InvIred}
Let $X$ be a compact metrizable space. Then we may write $X$ as the inverse limit of an inverse sequence $\mathbf{Q}=(\vert Q_i\vert,q_i^{i+1})$ of compact metric polyhedra, where each bonding map $q_i^{i+1}$ is $Q_i$-irreducible. \hfill $\square$
\end{lemma}

\begin{lemma}
\label{3.4[R-S]}Let $X$ be a compact metrizable space. Then there exists an inverse sequence $\mathbf{K}=(\vert K_i\vert, p_i^{i+1})$ of compact metric polyhedra $(\vert K_i\vert, d_i)$ along with a sequence of stability $(\ga_i)$ for $\mathbf{K}$ such that $\lim \mathbf{K}=X$, and for each $i\in \N$, $\mesh K_i<\ga_i$. We may also specify that for some $m \in \N$, whenever $i\geq m$, then $p_i^{i+1}:\vert K_{i+1}\vert \ra \vert K_i\vert$ is a $K_i$-irreducible simplicial map.
\end{lemma}

\noindent \textit{Proof}: Write $X=\lim \mathbf{Q}$, where $\mathbf{Q}=(\vert Q_i\vert,q_i^{i+1})$ is an inverse sequence of compact metric polyhedra $(\vert Q_i\vert,d_i)$ as in Lemma~\ref{InvIred}. By Lemma~\ref{stability}, we know that there is a sequence of stability $(\rho_i)$ for $\mathbf{Q}$. For each $i$, put $\ga_i=\rho_i/2$. Note that $(\ga_i)$ is also a sequence of stability for $\mathbf{Q}$.

Let $K_1$ be a subdivision of $Q_1$ with $\mesh K_1<\ga_1$. Suppose that $i\in \N$ and for each $1\leq j\leq i$, we have chosen a subdivision $K_j$ of $Q_j$ with $\mesh K_j<\ga_j$ and, when $1<j$, a map $p_{j-1}^{j}:\vert K_j\vert \ra \vert K_{j-1}\vert$ which is a simplicial approximation to $q_{j-1}^j$. Then select a subdivision $K_{i+1}$ of $Q_{i+1}$ with $\mesh K_{i+1}<\ga_{i+1}$, and which supports a simplicial approximation $p_i^{i+1}: \vert K_{i+1}\vert \ra \vert K_i\vert$ of $q_i^{i+1}$. Note that $d_i(q_i^{i+1},p_i^{i+1})<\ga_i$.

Let us check that $\mathbf{K}:=(\vert K_i\vert,p_i^{i+1})$ and $m=1$ satisfy all of the requirements.
Clearly $X=\lim \mathbf{K}$, since $(\ga_i)$ is a sequence of stability for $\mathbf{Q}$. It remains to show that the new bonding maps $p_i^{i+1}$ are $K_i$-irreducible. First, note that $q_i^{i+1}$ being $Q_i$-irreducible implies that $q_i^{i+1}$ is also $K_i$-irreducible.
Since $p_i^{i+1}$ is a simplicial approximation  of $q_i^{i+1}$, $p_i^{i+1}$ is a $K_i$-modification of $q_i^{i+1}$. By Lemma~\ref{K-irred},  $p_i^{i+1}$ is $K_i$-irreducible too. \hfill $\square$

\begin{defi}
\label{represent}
\emph{Whenever $X$ is a compact metrizable space, then we shall refer to an inverse sequence $\mathbf{K}$ of metric polyhedra $(\vert K_i\vert, d_i)$
which admits a sequence $(\ga_i)$ of positive numbers and $m \in \N$ so that the properties of Lemma~\ref{3.4[R-S]} are satisfied as a \emph{representation} of $X$ which is \emph{stable and simplicially irreducible from index $m$} with \emph{associated sequence of stability} $(\ga_i)$.}
\end{defi}

Of course, Lemma~\ref{3.4[R-S]} and its proof show that every  compact metrizable space $X$ has a representation $\mathbf{K}$ which is stable and simplicially irreducible from index $m=1$.

\vspace{2mm}

Next, we want to define a certain procedure which when applied to such $\mathbf{K}=\mathbf{K}_0$ as in Definition~\ref{represent} results in a $\mathbf{K}_1$ which is also a stable and simplicially irreducible (from some index $m$) representation of $X$. We will then show that if this procedure is repeated recursively in a controlled manner, resulting in a sequence $\mathbf{K}_1, \mathbf{K}_2, \ldots$, then there will be a limit $\mathbf{K}_\infty= \displaystyle{\lim_{j\ra \infty}} (\mathbf{K}_j)$ which also will be a representation of $X$.

\begin{lemma}
\label{ssirred}
Let $(\e_i)$ be a sequence of positive numbers. Let $X$ be a compact metrizable space, let $\mathbf{K}=(\vert K_i\vert,p_i^{i+1})$ be a representation of $X$ which is stable and simplicially irreducible from index $m_1$ with  an associated sequence of stability $(\ga_i)$, and
let $m \in \N_{\geq m_1}$. Define $\ga_i'=\ga_i$ if $1\leq i<m$, $\ga_m'=\frac{1}{2} [\ga_m - \mesh K_m ]$, and $\ga_i'=\ga_i/2$ if $i>m$.
Let $\Sa$ be a subdivision of $K_m$ with $\mesh \Sa <\min \{\e_m,\ga_m'\}$.
Then there exists an inverse sequence  $\mathbf{L}=(\vert L_i\vert, l_i^{i+1})$ as follows:
\begin{enumerate}
\item[(a)] in case $1\leq i < m$, then  $L_i=K_i$ and $l_i^{i+1}=p_i^{i+1}$,
\item[(b)] $L_m=\Sa$,
\item[(c)] for each $i\geq m+1$, $L_i$ is a subdivision of $K_i$ with $\mesh L_i< \min \{\e_i,\ga_i'\}$, and
\item[(d)] if $i\geq m+1$, $l_{i-1}^i : \vert L_i\vert \ra \vert L_{i-1}\vert$ is a simplicial approximation to the map $p_{i-1}^i$.\hfill $\square$
\end{enumerate}
\end{lemma}
\begin{defi}
\emph{
We shall call a pair $(\mathbf{L}, (\ga_i'))$ as in Lemma~\ref{ssirred} an $m$-\emph{shift} of $(\mathbf{K}, (\ga_i))$ from $\Sa$.}
\end{defi}

 Observe that $d_m(p_m^{m+1},l_m^{m+1})\leq \mesh \Sa < \frac{1}{2}[\ga_m - \mesh K_m]=\ga_m'$. Hence if $g:\vert L_{m+1}\vert \ra \vert L_m\vert$ is a map and $d_m(g,l_m^{m+1})<\ga_m'$, we may conclude that $d_m(g,p_m^{m+1})<\ga_m$. Indeed, the following is true:
\begin{enumerate}
\item[(e)] for each $i$, if $g:\vert L_{i+1}\vert \ra \vert L_i\vert$ is a map and $d_i(g,l_i^{i+1})<\ga_i'$, then $d_i(g,p_i^{i+1})<\ga_i$.
\end{enumerate}

Therefore we conclude:

\begin{lemma}
Whenever $(\mathbf{L}, (\ga_i'))$ is an $m$-\emph{shift} of $(\mathbf{K}, (\ga_i))$ from $\Sa$, then $\mathbf{L}$ is a stable and simplicially irreducible representation of $X$ from index $m$ with associated sequence of stability $(\ga_i')$.
\hfill $\square$
\end{lemma}

By exercising some additional care in the construction of $\mathbf{L}$, we may guarantee that
 for all $i$, $d_i(p_i^{i+1},l_i^{i+1})<\e_i$ (of course, $p_i^{i+1}=l_i^{i+1}$ if $i<m$).

\vspace{2mm}

It is routine to check that the next lemma holds true.

\begin{lemma}
\label{bigstate}
Let $X$ be a compact metrizable space, and let $\mathbf{K}_0$ be a representation of $X$ which is stable and simplicially irreducible from index $m_1$, with $(\ga_{(0),i})$ a sequence of stability. For every $m_1$-shift $(\mathbf{K}_1, (\ga_{(1),i}))$ of $(\mathbf{K}_0, (\ga_{(0),i}))$ from $\Sa_1$ (an appropriate subdivision of the triangulation of the $m_1$-term of $\mathbf{K}_0$), $\mathbf{K}_1$ is a representation of $X$ which is stable and simplicially irreducible from index $m_1$, with $(\ga_{(1),i})$ an associated sequence of stability. It satisfies property \emph{(e)}  with $(\ga_i')=(\ga_{(1),i})$ and $(\ga_i)=(\ga_{(0),i})$. The terms (as metric spaces) in $\mathbf{K}_0$ and $\mathbf{K}_1$ are equal. For $i< m_1$, $\ga_{(1),i}=\ga_{(0),i}$, the terms with index $i$ have the same triangulations in $\mathbf{K}_0$ and $\mathbf{K}_1$, and the bonding maps in $\mathbf{K}_0$ and $\mathbf{K}_1$ with subscript $i$ are equal. For $i \geq m_1$, $\ga_{(1),i}$ need not equal $\ga_{(0),i}$, the triangulation of the term in $\mathbf{K}_1$ with index $i$ is a subdivision of that in $\mathbf{K}_0$ with the same index, and the bonding map with subscript $i$ in $\mathbf{K}_1$ may differ from that in $\mathbf{K}_0$ with subscript $i$.

If $i_0\in \N$, $m_1<\ldots <m_{i_0}$ is a finite sequence in $\N$, and successively we have chosen $(\mathbf{K}_j, (\ga_{(j),i}))$ an $m_j$-shift of $(\mathbf{K}_{j-1}, (\ga_{(j-1),i}))$ from $\Sa_j$ (an appropriate subdivision of the $m_j$-term of $\mathbf{K}_{j-1}$), $1\leq j \leq i_0$, then we may conclude that $\mathbf{K}_{i_0}$ is a representation of $X$ which is stable and simplicially irreducible from index $m_{i_0}$, with $(\ga_{(i_0),i})$ an associated sequence of stability; it satisfies property \emph{(e)}  with $(\ga_i')=(\ga_{(i_0),i})$  and $(\ga_i)=(\ga_{(i_0-1),i})$. The terms (as metric spaces) in $\mathbf{K}_0$ and $\mathbf{K}_{i_0}$ are equal. For $i< m_{i_0}$, $\ga_{(i_0),i}=\ga_{(i_0-1),i}$, the terms with index $i$ have the same triangulations in $\mathbf{K}_{i_0-1}$ and $\mathbf{K}_{i_0}$, and the bonding maps in $\mathbf{K}_{i_0-1}$ and $\mathbf{K}_{i_0}$ with subscript $i$ are equal. For $i \geq m_{i_0}$, $\ga_{(i_0),i}$ need not equal $\ga_{(i_0-1),i}$, the triangulation of the term in $\mathbf{K}_{i_0}$ with index $i$ is a subdivision of that in $\mathbf{K}_{i_0-1}$ with the same index, and the bonding map with subscript $i$ in $\mathbf{K}_{i_0}$ may differ from that in $\mathbf{K}_{i_0-1}$ with subscript $i$.\hfill $\square$
\end{lemma}

Henceforth we typically shall write $(\vert K_{(j),i}\vert, p_{(j),i}^{i+1})$ to denote such a representation $\mathbf{K}_j$, $0\leq j \leq i_0$.
One should note that, whenever $i_0\geq j_0\geq j \geq 1$, then $K_{(j),m_j}=K_{(j_0),m_j}=\Sa_j$ when this occurs from the procedure in Lemma~\ref{bigstate}.

\begin{defi}
\label{Kinfty}
\emph{Let $X$ be a compact metrizable space and let $r:\N \ra \N$ be an increasing function. Let $\mathbf{K}_0$ be a representation of $X$ which is stable and simplicially irreducible from index $r(1)$, with $(\ga_{(0),i})$ a sequence of stability. Suppose that $(\mathbf{K}_j,(\ga_{(j),i}))$, $j \in \N$, is a sequence such that for each $j$, $(\mathbf{K}_j,(\ga_{(j),i}))$ is an $r(j)$-shift of $(\mathbf{K}_{j-1},(\ga_{(j-1),i}))$ from $\Sa_j$.}

\emph{Then for each $k\in \N$, if $m$, $l$, and $i$ are chosen so that $m\geq l \geq r(k)>i$, one sees that $p_{(l),i}^{i+1}=p_{(m),i}^{i+1}$ and $\ga_{(l),i}=\ga_{(m),i}$. So for each $i$, the sequences $(\ga_{(j),i})_{j \in \N}$ and $(p_{(j),i}^{i+1})_{j\in \N}$ are eventually constant.
Hence we may define an inverse sequence $\mathbf{K}_\infty=(\vert K_{(\infty),i}\vert, p_{(\infty),i}^{i+1})=\displaystyle{\lim_{j\ra \infty}}\mathbf{K}_j$ and a sequence $(\ga_{(\infty),i})=\displaystyle{\lim_{j\ra \infty}}(\ga_{(j),i})$ of positive numbers by putting $K_{(\infty),i}=\displaystyle{\lim_{j\ra \infty}} K_{(j),i}$ and $p_{(\infty),i}^{i+1}=\displaystyle{\lim_{j\ra \infty}}p_{(j),i}^{i+1}$.}
\end{defi}


From our construction and this definition, we can deduce the following:
\begin{lemma}
\label{KinftyX}
Assume the notation of Definition~\emph{\ref{Kinfty}}. Then $\mathbf{K}_\infty$ is a representation of $X$. If $i \in \N$, $g: \vert K_{(\infty),i+1}\vert \ra \vert K_{(\infty),i}\vert$ is a map, and $d_i(g,p_{(\infty),i}^{i+1})<\ga_{(\infty),i}$, then $d_i(g,p_{(0),i}^{i+1})<\ga_{(0),i}$  and hence
$(\ga_{(\infty),i})$ is a sequence of stability for $\mathbf{K}_\infty$.
\end{lemma}


\noindent \textit{Proof}: To show that $\mathbf{K}_\infty$ is a representation of $X$, it is enough to check that for all $i\in\N$, $d_i(p_{(\infty),i}^{i+1},p_{(0),i}^{i+1}) < \ga_{(0),i}$.

Take an $i\in \N$. If $i < r(1)$, then $p_{(\infty),i}^{i+1}=p_{(0),i}^{i+1}$ and $\ga_{(\infty),i}=\ga_{(0),i}$. Hence the statement $d_i(g,p_{(\infty),i}^{i+1})<\ga_{(\infty),i}$ implies that $d_i(g,p_{(0),i}^{i+1})<\ga_{(0),i}$.

\vspace{2mm}

If $i \geq r(1)$, then we know that $r(k-1)\leq i <r(k)$ for some $k\in\N_{\geq 2}$. The fact that $i <r(k)$ implies that $p_{(\infty),i}^{i+1}=p_{(k-1),i}^{i+1}$.
On the other hand, $r(k-1)\leq i$ implies that $\ga_{(j),i}$ has changed in every step of the construction from step $0$ to $(k-1)$. That is,
$\ga_{(j),i}\leq \frac{1}{2}\ga_{(j-1),i}$ for all $1\leq j\leq k-1$, so $\ga_{(j),i}\leq \frac{1}{2^j}\ga_{(0),i}$.
Therefore
\begin{align*}
d_i(p_{(\infty),i}^{i+1},p_{(0),i}^{i+1}) &=d_i(p_{(k-1),i}^{i+1},p_{(0),i}^{i+1}) \leq d_i(p_{(k-1),i}^{i+1},p_{(k-2),i}^{i+1})
+\ldots + d_i(p_{(1),i}^{i+1},p_{(0),i}^{i+1}) \\
&<\ga_{(k-1),i}+\ldots + \ga_{(1),i} \leq \frac{\ga_{(0),i}}{2^{k-1}} + \ldots + \frac{\ga_{(0),i}}{2}<\ga_{(0),i}\cdot \sum_{k=1}^\infty \frac{1}{2^k}=\ga_{(0),i}.
\end{align*}
By Lemma~\ref{stability}, $\lim \mathbf{K}_\infty =X$.

\vspace{2mm}

It remains to show that $d_i(g,p_{(\infty),i}^{i+1})<\ga_{(\infty),i}$ implies $d_i(g,p_{(0),i}^{i+1})<\ga_{(0),i}$. The fact that $i<r(k)$ implies that $\ga_{(\infty),i}=\ga_{(k-1),i}$. So
$d_i(g,p_{(\infty),i}^{i+1})=d_i(g,p_{(k-1),i}^{i+1})<\ga_{(k-1),i}$. Therefore
\begin{align*}
d_i(p_{(0),i}^{i+1},g) &\leq  d_i(p_{(0),i}^{i+1},p_{(1),i}^{i+1})+ d_i(p_{(1),i}^{i+1},p_{(2),i}^{i+1})+\ldots +d_i(p_{(k-2),i}^{i+1},p_{(k-1),i}^{i+1}) +d_i(p_{(k-1),i}^{i+1},g)\\
&< (\ga_{(1),i}+ \ga_{(2),i}+ \ldots + \ga_{(k-1),i}) + \ga_{(k-1),i}\\
&\leq \ga_{(0),i}\cdot \left(\left(\frac{1}{2} + \frac{1}{2^2} + \ldots + \frac{1}{2^{k-1}} \right) + \frac{1}{2^{k-1}}\right)=\ga_{(0),i}.
\end{align*}
\hfill $\square$

\section{Proof of the Main Theorem}

Let us now prove  Theorem~\ref{T}.



\vspace{2mm}

\noindent \textit{Proof}: We will construct, using induction:
\begin{enumerate}
\item[$\diamond$] an increasing function $r: \N \ra \N$,

\item[$\diamond$] sequences of numbers $(\da(i))_{i\in \N}$ and
$(\e(i))_{i\in \N}$ such that $0 < \e(i)<\frac{\da(i)}{3}<1$, for
all $i$,

\item[$\diamond$] a sequence of inverse sequences
$\mathbf{K}_j=(\vert K_{(j),i}\vert, p_{(j),i}^{i+1})$, for $j \in
\Z_{\geq0}$, as described in Lemma~\ref{bigstate}, with terms that are compact polyhedra and with surjective bonding maps, and with $\lim \mathbf{K}_j = X$
(in fact, these sequences are representations for $X$ that are stable and simplicially irreducible from index $r(j)$,
with stability sequences $(\ga_{(j),i})$,
and $\vert K_{(j),i}\vert=\vert K_{(0),i}\vert$, for all $i$ and
$j$ in $\N$),

\item[$\diamond$] a sequence of subdivisions $\Sa_i$ of
$K_{(i-1),r(i)}$, for $i \in \N$, and

\item[$\diamond$] a sequence of maps $g_{r(i-1)}^{r(i)}: \ \vert
\Sa_i^{(n+1)}\vert \ra \vert \Sa_{i-1}^{(n)}\vert$, for $i\geq 2$,
\end{enumerate}
such that for each $i$ for which the statement makes sense, we
have:
\begin{enumerate}
\item[(I)$_i$] $g_{r(i-1)}^{r(i)}$ and
$p_{(i-1),r(i-1)}^{r(i)}\vert_{\vert \Sa_i^{(n+1)}\vert} $ are
$\frac{\e(i-1)}{3}$ - close,

\item[(II)$_i$] for any $y \in \vert K_{(i-1),r(i)}\vert = \vert
\Sa_i \vert$, $\ \diam \ (p_{(i-1),r(i-1)}^{r(i)}(B_{\da(i)}(y)))
\ < \ \frac{\e(i-1)}{3}$,

\item[(III)$_i$] for $i>j$ and for any $y \in \vert
K_{(i-1),r(i)}\vert = \vert \Sa_i \vert$, $\ \diam \
(p_{(j),r(j)}^{r(i)}(B_{\e(i)}(y))) \ < \ \frac{\e(j)}{2^i}$,

\item [(IV)$_i$] $\mesh \Sa_i <  \min \
\{\frac{\e(i)}{3},\ga_{(i-1),r(i)}\}$, so $\mesh \Sa_i <  \e(i)$,
and

\item [(V)$_i$] for any $y \in \vert K_{(i-1),r(i)}\vert = \vert
\Sa_i \vert$, $\ B_{\e(i)}(y)\subset P_{y,i}\subset
B_{\da(i)}(y)$, where $P_{y,i}$ is a contractible subpolyhedron of
$\vert \Sa_i \vert$.
\end{enumerate}

\vspace{2mm}

\noindent In fact, this will prepare us to use Walsh's Lemma~\ref{specgenWa} with
$$\mathbf{X}=(\vert K_{(0),r(i)}\vert, \ p_{(i),r(i)}^{r(i+1)}), \
\ \mathbf{Z}=(\vert \Sa_i^{(n)}\vert, \
g_{r(i)}^{r(i+1)}\vert_{\vert \Sa_{i+1}^{(n)}\vert}).$$

\vspace{2mm}

Let us start the construction by taking a  representation for $X$ which is stable and simplicially
irreducible from index $1$: $\mathbf{K}_0 = (\vert
K_{(0),i}\vert, p_{(0),i}^{i+1})$, $\lim \mathbf{K}_0 = X$, with
stability sequence $(\ga_{(0),i})$.

Define $r(1):=1$.

We will choose $0< \da(1)<1$ any way we want. Next, we pick an
intermediate subdivision $\widetilde\Sa_1$ of $K_{(0),1}$ so that
for any $y \in \vert K_{(0),1}\vert$, any closed
$\widetilde\Sa_1$-vertex star containing $y$ is contained in the
closed $\da(1)$-ball $B_{\da(1)}(y)$. (A closed
$\widetilde\Sa_1$-vertex star is a closed star $\overline{\st}
(w,\widetilde\Sa_1)$ in the complex $\widetilde\Sa_1$ whose center $w$ is a vertex of $\widetilde\Sa_1$.) It is enough to make $\mesh \widetilde\Sa_1 < \frac{\da(1)}{2}$, so $\diam (\overline{\st}
(w,\widetilde\Sa_1)) \ \leq \ 2 \mesh \widetilde\Sa_1 \ < \da(1)$).

Now choose an $\e(1)$ so that $0<\e(1)<\frac{\da(1)}{3}$, and  for
any $y \in \vert K_{(0),1}\vert$, the closed $\e(1)$-ball
$B_{\e(1)}(y)$ sits inside an open vertex star with respect to
$\widetilde\Sa_1$. (This can be done as follows: form the open
cover for $\vert K_{(0),1}\vert$ consisting of the open stars
$\st(w,\widetilde\Sa_1)$. There is a Lebesgue number $\la$ for
this cover, so make your $\e(1) < \frac{\la}{2}$. Then for any $y
\in \vert K_{(0),1}\vert$, $\diam B_{\e(1)}(y)<\la \ \Rightarrow \
B_{\e(1)}(y)\subset \st (w_0,\widetilde\Sa_1)$, for some $w_0 \in
\widetilde\Sa_1^{(0)}$. Fix such $w_0$ for each $y$.)

Note that for any $y \in \vert K_{(0),1}\vert$,
$B_{\e(1)}(y)\subset \ \vert \overline{\st}
(w_0,\widetilde\Sa_1)\vert \ \subset B_{\da(1)}(y)$. Define
$P_{y,1} := \vert \overline{\st} (w_0,\widetilde\Sa_1)\vert$,
which is a contractible subpolyhedron of $\vert K_{(0),1}\vert$,
so (V)$_1$ is satisfied.

Choose a subdivision $\Sa_1$ of $\widetilde\Sa_1$  with $\mesh
\Sa_1 <\min \ \{ \frac{\e(1)}{3}, \ga_{(0),1}\}$, which implies
(IV)$_1$.

Let $(\mathbf{K}_1,(\ga_{(1),i}))$ be a $1$-shift of
$(\mathbf{K}_0,(\ga_{(0),i}))$ from $\Sa_1$, i.e.,
$\mathbf{K}_1=(\vert K_{(1),i}\vert, p_{(1),i}^{i+1})$ is an
inverse sequence with $K_{(1),1}=\Sa_1$, limit equal $X$, and
stability sequence $(\ga_{(1),i})$.
Note that at this point, all bonding maps in $\mathbf{K}_1$ are simplicial because $\mathbf{K}_1$ is simplicially irreducible from index $1$.
This concludes the basis of induction.


\vspace{1mm}

\textit{Step of induction.} Let $k\in \N_{\geq 2}$. Suppose that we have chosen, as
required above,
\begin{enumerate}
\item[$\diamond$] for $j=1,\dots , k-1$, the numbers $r(j)$,
$\da(j)$, $\e(j)$,

\item[$\diamond$] for $j=0,\dots , k-1$, the inverse sequences
$\mathbf{K}_j=(\vert K_{(j),i}\vert, p_{(j),i}^{i+1})$, which are stable and simplicially irreducible from index $r(j)$, with
stability sequences $(\ga_{(j),i})$,

\item[$\diamond$] for $j=1,\dots , k-1$, subdivisions $\Sa_j$ of
$K_{(j-1),r(j)}$, and

\item[$\diamond$] for $j=2,\dots , k-1$, maps $g_{r(j-1)}^{r(j)}:
\ \vert \Sa_j^{(n+1)}\vert \ra \vert \Sa_{j-1}^{(n)}\vert$,
\end{enumerate}
so that the properties (I)$_j$-(V)$_j$ are satisfied for each
$j=1,\dots , k-1$ for which they make sense.

Focus on the inverse sequence $\mathbf{K}_{k-1}=(\vert
K_{(k-1),i}\vert, \ p_{(k-1),i}^{i+1})$. For $i\geq r(k-1)$, the bonding maps $p_{(k-1),i}^{i+1}$ are simplicial. Recall that $\lim
\mathbf{K}_{k-1}=X$, and notice that $K_{(k-1),r(k-1)}=\Sa_{k-1}$.
Let

\vspace{-5mm}

$$\mathbf{Y}_{k-1}:=(\vert K_{(k-1),i}^{(n+1)}\vert, \
p_{(k-1),i}^{i+1}\vert_{\ \vert K_{(k-1),i+1}^{(n+1)}\vert
})_{_{i\geq r(k-1)}}$$


\noindent be the inverse sequence of the
$(n+1)$-skeleta of the polyhedra in $\mathbf{K}_{k-1}$, starting
with the $(r(k-1))$-th polyhedron onward, where the
bonding maps are the restrictions of the original bonding maps. Notice
that every $p_{(k-1),i}^{i+1}\vert_{\ \vert K_{(k-1),i+1}^{(n+1)}\vert
}:\vert K_{(k-1),i+1}^{(n+1)}\vert\ra \vert K_{(k-1),i}^{(n+1)}\vert $
is still simplicial and surjective: since $p_{(k-1),i}^{i+1}$ is
simplicial and surjective, for every simplex $\sa \in
K_{(k-1),i}^{(n+1)}$ with $\dim \sa = k$, there exists a simplex
$\tau \in K_{(k-1),i+1}$ such that $\dim \tau\geq k$ and
$p_{(k-1),i}^{i+1}\ (\tau) = \sa$. So there must be a $k$-face of
$\tau$ which is mapped by $p_{(k-1),i}^{i+1}$ onto $\sa$. In
particular, for every $(n+1)$-dimensional $\sa \in
K_{(k-1),i}^{(n+1)}$, there exists an $(n+1)$-simplex in
$K_{(k-1),i+1}$ that is mapped onto $\sa$ by $p_{(k-1),i}^{i+1}$.

Now let $Y_{k-1}=\lim \mathbf{Y}_{k-1}$. Then $\dim Y_{k-1}\leq
n+1$, because $\dim \vert K_{(k-1),i}^{(n+1)}\vert\leq n+1$, and $X\tau K$ implies $Y_{k-1} \ \tau K$, because $Y_{k-1}\subset X$. So by
Lemma~\ref{dimGY}, we get $\dim_G Y_{k-1} \leq n$. Since $P_G=\P$, Lemma~\ref{PGPrimesZ} implies
$\dim_\Z Y_{k-1} = \dim_G Y_{k-1} \leq n$, so we can apply
Edwards' Theorem~\ref{Ed} to $\mathbf{Y}_{k-1}$, noticing that the first
entry in $\mathbf{Y}_{k-1}$ has index $r(k-1)$.

So there exists an $s\in \N$, $s>r(k-1)$ and a map
$\widehat{g}_{r(k-1)}^s : \vert K_{(k-1),s}^{(n+1)}\vert \ra \vert
K_{(k-1),r(k-1)}^{(n)}\vert$ so that if $z \in \vert
K_{(k-1),s}^{(n+1)}\vert$, and $p_{(k-1),r(k-1)}^s(z)$ lands in
the combinatorial interior $\overset{\circ}\sa$ of a simplex $\sa$
of $K_{(k-1),r(k-1)}^{(n+1)}$, then $\widehat{g}_{r(k-1)}^s (z)$
lands in $\sa$. This will help us get the property (I)$_k$.

\begin{displaymath}
\xymatrix{
\vert K_{(k-1),r(k-1)}^{(n)}\vert \ar@{_{(}->}[d] & & &                                 &               &\\
 \vert K_{(k-1),r(k-1)}^{(n+1)}\vert & & & \vert K_{(k-1),r(k)}^{(n+1)}\vert \ar[lll]^{p_{(k-1),r(k-1)}^{r(k)}\vert }\ar@{-->}[lllu]_{\widehat{g}_{r(k-1)}^{r(k)}} & \cdots \ar[l] & Y_{k-1}
}
\end{displaymath}

Define $r(k):= s$. Using the uniform continuity of the map
$p_{(k-1),r(k-1)}^{r(k)}$, choose $0< \da(k)< 1$ so that (II)$_k$
is true: $$\forall y \in \vert K_{(k-1),r(k)}\vert ,\ \ \diam \ (
p_{(k-1),r(k-1)}^{r(k)}(B_{\da(k)}(y))) \ < \ \frac{\e(k-1)}{3}.$$
Pick an intermediate subdivision $\widetilde\Sa_k$ of
$K_{(k-1),r(k)}$ so that for any $y \in \vert
K_{(k-1),r(k)}\vert$, any closed $\widetilde\Sa_k$-vertex star
containing $y$ is contained in $B_{\da(k)}(y)$.

Now choose an $\e(k)$ so that $0<\e(k)<\frac{\da(k)}{3}$, and so
that (III)$_k$ and (V)$_k$ will hold, namely:\\ first make sure that for
all $y \in \vert K_{(k-1),r(k)}\vert$, the closed $\e(k)$-ball
centered at $y$ sits inside an open $\widetilde\Sa_k$-vertex star,
i.e., $B_{\e(k)}(y) \subset \st (w_0,\widetilde\Sa_k)$, for some
$w_0 \in \widetilde\Sa_k^{(0)}$. Therefore $B_{\e(k)}(y)\subset \
\vert \overline{\st} (w_0,\widetilde\Sa_k)\vert \ \subset
B_{\da(k)}(y)$. Define $P_{y,k} := \vert \overline{\st}
(w_0,\widetilde\Sa_k)\vert$, which is a contractible subpolyhedron
of $\vert K_{(k-1),r(k)}\vert$, so (V)$_k$ is satisfied. Next, we
know that for all $j<k$, the maps $p_{(j),r(j)}^{r(k)}$ are
uniformly continuous. We also know that, in our notation, $j<k$
implies that $p_{(j),r(j)}^{r(k)}=p_{(k-1),r(j)}^{r(k)}$. So we
can make a choice of $\e(k)$ so that we have: for any $y \in \vert
K_{(k-1),r(k)}\vert$,
\begin{align*}
\diam \ ( p_{(1),r(1)}^{r(k)}(B_{\e(k)}(y))) \ &< \
\frac{\e(1)}{2^k},\\
\diam \ ( p_{(2),r(2)}^{r(k)}(B_{\e(k)}(y))) \ &< \
\frac{\e(2)}{2^k},\\
\vdots\\
\diam \ ( p_{(k-1),r(k-1)}^{r(k)}(B_{\e(k)}(y))) \ &< \
\frac{\e(k-1)}{2^k}.
\end{align*}
So (III)$_k$ is true.

Choose a subdivision $\Sa_k$ of $\widetilde\Sa_k$  with $\mesh
\Sa_k < \ga_{(k-1),r(k)}$, where $\ga_{(k-1),r(k)}$ is from the
stability sequence $(\ga_{(k-1),i})$ for $\mathbf{K}_{k-1}$. Also
make sure that $\mesh \Sa_k <\frac{\e(k)}{3}$, which implies
(IV)$_k$. Note that $\Sa_k$ is a subdivision of $K_{(k-1),r(k)}$.
\begin{displaymath}
\xymatrix{
\mathbf{K}_{k-1}: & & \cdots &\vert K_{(k-1),r(k)}\vert  \ar[l] && \cdots\hspace{1.5cm}  \ar[ll]_{{p_{(k-1),r(k)}^{r(k)+1}} } & &\\
\mathbf{K}_{k}:& \cdots &  \vert \Sa_k\vert=\hspace{-15mm}\ar[l] &\vert K_{(k),r(k)}\vert  \ar[u]^{id} \ar@/_/@{.>}[u]|{j}      &&  \vert K_{(k),r(k)+1}\vert  \ar[ll]_{{p_{(k),r(k)}^{r(k)+1}} }     & \cdots \ar[l]  &X\\
\mathbf{Y}_{k}: & &\vert \Sa_k^{(n+1)}\vert=\hspace{-9mm}& \vert K_{(k),r(k)}^{(n+1)}\vert \ar@{_{(}->}[u]                  &&  \vert K_{(k),r(k)+1}^{(n+1)}\vert  \ar@{_{(}->}[u]  \ar[ll]_{{p_{(k),r(k)}^{r(k)+1}}\vert }     & \cdots \ar[l]  & Y_k^{\  } \ar@{_{(}->}[u]\\
}
\end{displaymath}

Now we can build $\mathbf{K}_{k}=(\vert K_{(k),i}\vert, \
p_{(k),i}^{i+1})$ as an $r(k)$-shift of
$(\mathbf{K}_{k-1},(\ga_{(k-1),i}))$ from $\Sa_k$, i.e.,
$\mathbf{K}_k=(\vert K_{(k),i}\vert, p_{(k),i}^{i+1})$ is an
inverse sequence with $K_{(k),r(k)}=\Sa_k$ and limit $X$, and
stability sequence $(\ga_{(k),i})$. For index $i\geq r(k)$, the bonding maps $p_{(k),i}^{i+1}$ are simplicial.

Let $j:\vert \Sa_k\vert \ra \vert K_{(k-1),r(k)}\vert$ be a
simplicial approximation to the identity map. Since $j$ is
simplicial, $j(\vert \Sa_k^{(n+1)}\vert )\subset \vert
K_{(k-1),r(k)}^{(n+1)}\vert$, so treat
$j\vert_{\vert \Sa_k^{(n+1)}\vert}:\vert \Sa_k^{(n+1)}\vert \ra \vert K_{(k-1),r(k)}^{(n+1)}\vert$.

Define $g_{r(k-1)}^{r(k)}:=\widehat{g}_{r(k-1)}^{r(k)}\circ
j\vert_{\vert \Sa_k^{(n+1)}\vert}: \vert \Sa_k^{(n+1)}\vert \ra
\vert K_{(k-1),r(k-1)}^{(n)}\vert= \vert \Sa_{k-1}^{(n)}\vert$.
For any $y \in \vert \Sa_k^{(n+1)}\vert$, $y$ and $j(y)$ have to
be contained in the same simplex of $K_{(k-1),r(k)}$. Since
$p_{(k-1),r(k-1)}^{r(k)}: \vert K_{(k-1),r(k)}\vert \ra \vert K_{(k-1),r(k-1)}\vert$ is simplicial,
$p_{(k-1),r(k-1)}^{r(k)}(y)$ and $p_{(k-1),r(k-1)}^{r(k)}(j(y))$
land in the same simplex $\tau$ of $K_{(k-1),r(k-1)}=\Sa_{k-1}$.
On the other hand, because of our choice of
$\widehat{g}_{r(k-1)}^{r(k)}$, if $p_{(k-1),r(k-1)}^{r(k)}(j(y))$
lands in $\overset{\circ}\sa$, for some simplex $\sa$ of
$K_{(k-1),r(k-1)}^{(n+1)}$ which is a face of $\tau$, then
$\widehat{g}_{r(k-1)}^{r(k)} (j(y))$ lands in $\sa$, too.
Therefore
$$d_{k-1}(p_{(k-1),r(k-1)}^{r(k)}(y),\widehat{g}_{r(k-1)}^{r(k)} (j(y))) \
\leq \ \mesh K_{(k-1),r(k-1)}=\mesh \Sa_{k-1} <
\frac{\e(k-1)}{3}.$$
Hence $g_{r(k-1)}^{r(k)}$ and
$p_{(k-1),r(k-1)}^{r(k)}\vert_{\vert \Sa_k^{(n+1)}\vert}$ are
$\frac{\e(k-1)}{3}$-close, so (I)$_k$ is true. This concludes the inductive step. The following diagram summarizes the preceding construction.

\begin{displaymath}
\xymatrix{
\vert \Sa_{k-1}^{(n)}\vert = \hspace{-13mm} & \vert K_{(k-1),r(k-1)}^{(n)}\vert \ar@{_{(}->}[d] & &            &&&\\
\vert \Sa_{k-1}^{(n+1)}\vert = \hspace{-9mm} & \vert K_{(k-1),r(k-1)}^{(n+1)}\vert \ar@{_{(}->}[d]& & \vert K_{(k-1),r(k)}^{(n+1)}\vert \ar@{_{(}->}[d] \ar[ll]^{p_{(k-1),r(k-1)}^{r(k)}\vert} \ar[llu]_{\widehat{g}_{r(k-1)}^{r(k)}} &&&\\
\vert \Sa_{k-1}\vert = \hspace{-14mm} & \vert K_{(k-1),r(k-1)}\vert & & \vert K_{(k-1),r(k)}\vert\ar[ll]^{p_{(k-1),r(k-1)}^{r(k)}} && \vert K_{(k),r(k)}^{(n+1)}\vert \ar[llu]_{j\vert} \ar@{_{(}->}[d] & \hspace{-9mm} =\vert \Sa_k^{(n+1)}\vert\\
 & && & & \vert K_{(k),r(k)}\vert \ar@/^/[llu]^{id} \ar[llu]_{j}& \\
}
\end{displaymath}

Notice that the inverse sequence $$\mathbf{X}:=(\vert
K_{(0),r(i)}\vert, \ p_{(i),r(i)}^{r(i+1)})\ =\ (\vert
K_{(i),r(i)}\vert, \ p_{(i),r(i)}^{r(i+1)})\ =\ (\vert \Sa_i\vert,
\ p_{(i),r(i)}^{r(i+1)})\  $$ is a subsequence of
$\mathbf{K}_\infty = (\vert K_{(\infty),i}\vert, \
p_{(\infty),i}^{i+1})=(\vert K_{(0),i}\vert, \
p_{(\infty),i}^{i+1})$. By Lemma~\ref{KinftyX}, $\lim \mathbf{K}_\infty =X$, so $\lim \mathbf{X}\ $ is homeomorphic to
$X$. Without loss of generality, assume that $\lim \mathbf{X}=X$.

Let $\mathbf{Z}:=(\vert \Sa_i^{(n)}\vert, \
g_{r(i)}^{r(i+1)}\vert_{\vert \Sa_{i+1}^{(n)}\vert})$. Since
$\vert \Sa_i^{(n)}\vert$ are metrizable, compact and nonempty,
$\lim \mathbf{Z}=Z$ is a nonempty compact metrizable space.
Clearly, $\dim Z \leq n$, which also implies that $\dim_G Z \leq
n$. Now  $Z\tau K$ follows from Lemma~\ref{dimGY}.

Apply Walsh's Lemma~\ref{specgenWa} to these $\mathbf{X}$ and $\mathbf{Z}$: since the
requirements (I)-(VI) of Lemma~\ref{specgenWa} are satisfied, there is a
cell-like surjective map $\pi: Z \ra X$.\hfill $\square$

\vspace{2mm}

\begin{cor}
Let $G$ be an abelian group with $P_G=\P$. Let $K$ be a connected \emph{CW}-complex with
$\pi_1(K)\cong G$. Then every compact metrizable space $X$ with $X\tau K$ has to have $\dim X \leq 1$.
\end{cor}

\noindent\textit{Proof}: Theorem~\ref{T} is true for $n=1$, so for any compact metrizable space $X$ with $X\tau K$, we can find a compact metrizable space $Z$ with $\dim Z \leq 1$, $Z\tau K$ and a cell-like map $\pi:Z \ra X$. Note that cell-like maps are always surjective. Also, cell-like maps are $G$-acyclic, so in particular, $\pi$ is a $\Z$-acyclic map.

The Vietoris-Begle Theorem implies that a $G$-acyclic map cannot raise  $\dim_G$-dimension.
Since $\dim Z \leq 1$ implies that $\dim_\Z Z\leq 1$, and since $\pi$ is a $\Z$-acyclic map, we have that $\dim_\Z X \leq 1$, too. Recall that $\dim_\Z X \leq 1 \ \Leftrightarrow \ \dim X\leq 1$.\hfill $\square$


\begin{thebibliography}{99}

\bibitem[AJR]{AJR} S.~Ageev, R.~Jim\'{e}nez and L.~Rubin, \emph{Cell-like resolutions in the strongly countable $\Z$-dimensional case}, Topology and its Appls.\ \textbf{140} (2004), 5-14.

\bibitem[ARS]{ARS} S.~Ageev, D.~Repov\v{s}, E.~Shchepin, \emph{On the softness of the Dranishnikov resolution}, Proc.\ Steklov Inst.\ Math.\ \textbf{212} (1996), 3-27.

\bibitem[Al]{Al} P.\ S.~Aleksandrov, \emph{Einige Problemstellungen in der mengentheoretischen Topologie}, Mat.\ Sb.\ \textbf{43} (1936), 619-634.

\bibitem[Bi]{Bi} R.\ H.~Bing, \emph{The cartesian product of a certain nonmanifold and a line is $E^4$}, Ann.\ of Math.\ \textbf{70} (1959), 399-412.

\bibitem[Br]{Br} M.~Brown, \emph{Some applications of an approximation theorem for inverse limits}, Proc.\ Amer.\ Math.\ Soc.\ \textbf{11} (1960), 478-483.

\bibitem[Da]{Da} R.~Daverman, \emph{Detecting the disjoint disks property}, Pacific J.\ Math.\ \textbf{93} (1981), 277-298.

\bibitem[Dr1]{Dr1} A.~Dranishnikov, \emph{On P.S.~Aleksandrov's problem}, Mat.\ Sb.\ \textbf{135} (\textbf{4}) (1988), 551-557.

\bibitem[Dr2]{Dr2} A.~Dranishnikov, \emph{On homological dimension modulo $p$}, Math.\ USSR Sb.\ \textbf{60} \textbf{(2)} (1988), 413-425.

\bibitem[Dr3]{Dr3} A.~Dranishnikov, \emph{Cohomological dimension theory of compact metric spaces}, Topology Atlas Invited Contributions,
\texttt{http://at.yorku.ca/t/a/i/c/43.pdf}

\bibitem[Dr4]{Dr4} A.~Dranishnikov, \emph{Rational homology manifolds and rational resolutions}, Topology and its Appls.\ \textbf{94} (1999), 75-86.

\bibitem[Du]{Du} J.~Dugundji, \emph{Topology}, Allyn and Bacon, Boston, 1966.

\bibitem[Dy1]{Dy1} J.~Dydak, \emph{Cohomological dimension and metrizable spaces}, Trans.\ Amer.\ Math.\ Soc.\ \textbf{337} (1993), 219-234.

\bibitem[DW]{DW} J.~Dydak and J.~Walsh, \emph{Complexes that arise in cohomological dimension theory: a unified approach}, J.\ London Math.\ Soc.\ \textbf{(2) (48), no.\ 2} (1993), 329-347.

\bibitem[Dy]{Dy} E.~Dyer, \emph{On the dimension of products}, Fund.\ Math.\ \textbf{47} (1959), 141-160.

\bibitem[Ed]{Ed} R.\ D.~Edwards, \emph{A theorem and a question related to cohomological dimension and cell-like maps}, Notices Amer.\ Math.\ Soc.\ \textbf{25} (1978), A-259.

\bibitem[En]{En} R.~Engelking, \emph{General Topology}, Heldermann Verlag, Berlin, 1989.

\bibitem[Fr]{Fr} H.~Freudenthal, \emph{Entwicklungen von R\"{a}umen und ihren Gruppen}, Compositio Math.\ \textbf{4} (1937), 145-234.

\bibitem[Ha]{Ha} A.~Hatcher, \emph{Algebraic Topology}, Cambridge University Press, 2002.

\bibitem[HW]{HW} W.~Hurewicz and H.~Wallman, \emph{Dimension Theory}, Princeton University Press, 1948.

\bibitem[IR]{IR} I.~Ivan\v{s}i\'{c} and L.~Rubin, \emph{The extension dimension of universal spaces}, Glas.\ Mat.\ \textbf{38} (\textbf{58}) (2003), 121-127.

\bibitem[JR]{JR} R.~Jim\'{e}nez and L.~Rubin, \emph{An addition theorem for $n$-fundamental dimension in metric compacta}, Topology and its Appls.\ \textbf{62} (1995), 281-297.

\bibitem[Ku]{Ku} V.\ I.~Kuz'minov, \emph{Homological dimension theory}, Russian Math.\ Surveys \textbf{23} (1968), 1-45.

\bibitem[KY1]{KY1} A.~Koyama and K.~Yokoi, \emph{A unified approach of characterizations and resolutions for cohomological dimension modulo p}, Tsukuba J.\ Math.\ \textbf{18 (2)} (1994), 247-282.

\bibitem[KY2]{KY2} A.~Koyama and K.~Yokoi, \emph{Cohomological dimension and acyclic resolutions}, Topology and its Appls.\ \textbf{120} (2002), 175-204.

\bibitem[Le1]{Le1} M.~Levin, \textit{Acyclic resolutions for arbitrary groups}, Isr.\ J.\ Math.\ \textbf{135} (2003), 193-204.

\bibitem[Le2]{Le2} M.~Levin, \textit{Rational acyclic resolutions}, Algebraic and Geometric Topology \textbf{5} (2005), 219-235.

\bibitem[MR]{MR} S.~Marde\v{s}i\'{c} and L.~Rubin, \emph{Cell-like mappings and nonmetrizable compacta of finite cohomological dimension}, Trans.\ Amer.\ Math.\ Soc.\ \textbf{313} (1989), 53-79.

\bibitem[MS1]{MS1} S.~Marde\v{s}i\'{c} and J.~Segal, \emph{Shape theory}, North-Holland, Amsterdam, 1982.

\bibitem[MS2]{MS2} S.~Marde\v{s}i\'{c} and J.~Segal, \emph{Stability of almost commutative inverse systems of compacta}, Topology and its Appls.\ \textbf{31} (1989), 285-299.

\bibitem[Mu]{Mu} J.~Munkres, \emph{Topology}, Prentice Hall, Upper Saddle River, New Jersey, 1975.

\bibitem[RS1]{RS1} L.~Rubin and P.~Schapiro, \emph{Cell-like maps onto non-compact spaces of finite cohomological dimension}, Topology and its Appls.\ \textbf{27} (1987), 221-244.

\bibitem[RS2]{RS2} L.~Rubin and P.~Schapiro, \emph{Resolutions for metrizable compacta in extension theory}, Trans.\ Amer.\ Math.\ Soc.\ \textbf{358} (2005), 2507-2536.


\bibitem[Sp]{Sp} E.~Spanier, \emph{Algebraic Topology}, McGraw-Hill, New York, 1966.

\bibitem[Va]{Va} V.\ A.~Vassiliev, \emph{Introduction to topology}, Student Mathematical Library, V. 14, AMS 2001.

\bibitem[Wa]{Wa} J.~Walsh, \emph{Dimension, cohomological dimension, and cell-like mappings}, Shape Theory and Geometric Topology, Lecture Notes in Mathematics, volume 870, Springer Verlag, Berlin, 1981, pp.\ 105-118.

\end{thebibliography}
\end{document}